\documentclass[a4paper,11pt,reqno]{article}

\usepackage{amsmath,amsfonts,amsthm,amssymb}
\usepackage{mathpazo}
 
\newcommand{\ol}[1]{\overline{#1}}

\newcommand{\cB}{\mathcal{B}}

\newcommand{\cZ}{\mathcal{Z}}

\newcommand{\ott}{[0,T]}

\newcommand{\II}{\mathcal{I}}
\DeclareMathOperator{\id}{\text{Id}}

\newcommand{\eps}{\varepsilon}

\makeatletter
\newcommand{\eqcolon}{\mathrel{\mathord{=}\raise.2\p@\hbox{:}}}
\newcommand{\coloneq}{\mathrel{\raise.2\p@\hbox{:}\mathord{=}}}
\makeatother
\newcommand{\der}{\delta}

\newcommand{\RR}{\mathbb{R}}
\newcommand{\NN}{\mathbb{N}}

\newcommand{\FF}{\mathcal{F}}
\newcommand{\QQ}{\mathcal{Q}}

\newcommand{\CC}{\mathcal{C}}

\newcommand{\DD}{\mathcal{D}}
\newcommand{\ZZ}{\mathcal{Z}}
\newcommand{\LL}{\mathcal{L}}

\newcommand{\norm}[1]{\lVert #1\rVert}

\newcommand{\R}{\mathbb R}

\newcommand{\bTT}{\mathbb T}

\newcommand{\cac}{\mathcal C}

\newcommand{\cj}{\mathcal J}
\newcommand{\cl}{\mathcal L}

\newcommand{\lcl}{\left\{}
\newcommand{\rcl}{\right\}}

\newtheorem{theorem}{Theorem}[section]

\newtheorem{corollary}[theorem]{Corollary}

\newtheorem{definition}[theorem]{Definition}
\newtheorem{example}[theorem]{Example}

\newtheorem{lemma}[theorem]{Lemma}

\newtheorem{proposition}[theorem]{Proposition}

\newtheorem{remark}[theorem]{Remark}

\newcommand{\TT}{\mathcal{T}}
\newcommand{\cA}{\mathcal{A}}

\usepackage{pstricks} 
\usepackage{pst-node}
\usepackage{pst-tree}

\newcommand{\tnode}{\TR{\raisebox{0.5pt}{\ensuremath{\bullet}}}} 
\newcommand{\tnodel}[1]{\TR{\makebox[6pt][l]{\raisebox{0.5pt}{\ensuremath{\bullet{\text{\footnotesize
            #1}}}}}}}

\newcommand{\troot}{\bullet}

\newcommand{\tsnode}{\TR{\raisebox{0.5pt}{\ensuremath{\scriptstyle\bullet}}}} 
\newcommand{\tsroot}{\ensuremath{\scriptstyle\bullet}}
\newcommand{\aaabbb}{\pstree{\tsnode}{\pstree{\tsnode}{\tsnode}}} 
\newcommand{\aaaabbbb}{\pstree{\tsnode}{\pstree{\tsnode}{\pstree{\tsnode}{\tsnode}}}} 
\newcommand{\aabb}{\pstree{\tsnode}{\tsnode}} 
\newcommand{\aababb}{\pstree{\tsnode}{\tsnode \tsnode}}
\newcommand{\aaababbb}{\pstree{\tsnode}{\aababb}}
\newcommand{\aabababb}{\pstree{\tsnode}{\tsnode\tsnode\tsnode}}
\newcommand{\aaabbabb}{\pstree{\tsnode}{\tsnode \aabb}} 
\newcommand{\smalltrees}{\psset{levelsep=-5pt,nodesep=-2pt,treesep=1pt}}
\newcommand{\largetrees}{\psset{levelsep=-10pt,nodesep=-4pt,treesep=5pt}
}
\largetrees

\begin{document}

\title{Ramification of rough paths}
\author{
  { Massimiliano Gubinelli}              \\
{\small\it Equipe de probabilit\'es, statistique et mod\'elisation }\\[-0.1cm]
  {\small\it  Universit\'e Paris Sud}          \\[-0.1cm]
{\small \it 91 405 Orsay cedex, France}\\[-0.1cm] 
 {\small  {\tt gubinell@math.u-psud.fr}}   
}
\date{}
\maketitle

\begin{abstract}
The stack of iterated integrals of a path 
 is embedded in a larger algebraic structure where
iterated integrals are indexed by decorated rooted trees and where an
extended Chen's multiplicative property involves the
D\"urr-Connes-Kreimer coproduct on rooted trees.
This turns out to be the natural setting for a non-geometric theory of
rough paths. \\[0.1cm]

\end{abstract}

\section{Introduction}

Since the seminal work of Butcher on integration
methods~\cite{MR1993957,MR0305608} rooted trees (otherwise called
Cayley trees~\cite{Cayley:1881fk}) are recognized as a basic
combinatorial structure underlying the numerical and exact solution of
ordinary differential equations~(see for example~\cite{MR0403225}  and the monograph of
Hairer-N\o rsett-Wanner~\cite{MR1227985,MR0403225}). Trees are also
present in the work of Connes--Kreimer~\cite{MR1660199,MR1748177,MR1810779} on the combinatorial
structure of renormalization in perturbative Quantum Field Theory  and connections
between Runge-Kutta methods and renormalization has been explored by
Brouder~\cite{MR2106008,brouder1}. Connes and Kreimer explored a
Hopf algebra 
structure on rooted trees  to disentangle nested sub-divergences in the Feynman diagrams of perturbative QFT. Starting point is the work of Kreimer~\cite{MR1633004,MR1797019} which introduced nested integrals indexed by trees in the analysis of Feynman diagrams. The same Hopf algebra was described before by D\"ur~\cite{MR857100} (for basic results on Hopf algebras see
e.g.~\cite{MR0252485}).

Literature on combinatorial and algebraic properties of rooted trees is quite large, we prefer to single out the
work of  Hoffman~\cite{MR1990174} and the two papers of
Foissy~\cite{MR1909461,MR1905177} on labeled rooted trees.

A sub-algebra of the  Hopf algebra of rooted trees is isomorphic to the Hopf algebra of
Chen's iterated integrals~\cite{MR0454968,MR1847673} which is at the
base of Lyons theory of rough paths~\cite{MR1654527}. Lyons
theory allows to define and solve differential
equations driven by irregular ``noises''. For an 
exposition see the work of Lyons cited above, the book of Lyons and Qian~\cite{MR2036784}, the introductory article of Lejay~\cite{MR2053040}. For alternative approaches to rough paths see the paper~\cite{MR2091358} of the present author or Feyel-de~La~Pradelle~\cite{feyel}. 

 Chen~\cite{MR0454968}  showed that a given path in a
manifold can be encoded  in the Hopf algebra of its
iterated integrals. Lyons~\cite{MR1654527} realized that this encoding is good enough
to recover solutions of differential equation driven by such a path.

The aim of the present paper is to build a bridge between rooted trees
and rough paths. Here we would like to describe how to encode a
control path in a function on labeled rooted trees which we call a
\emph{branched rough path} and then generalize the theory of Lyons
to build solutions of driven differential equation by using this
new encoding.

The advantage of this approach is that we can dispose of
the notion of \emph{geometric} rough path which is fundamental in Lyons
theory. Geometric rough paths possess a rich structure and present
nice connections with the geometry of certain Carnot
groups~\cite{MR2144230} but there are situations where the geometric
property is not natural, e.g. in It\^o stochastic integration or in
infinite-dimensional generalizations of rough
paths~\cite{kdv,TindelGubinelli}. A more abstract motivation is to
prove that it is possible to build a complete theory of rough paths
(at any level of roughness) in the non-geometric setting. 
Series over trees can be helpful also in the geometric setting:
recently Neuenkirch--Nourdin--R\"o{\ss}ler--Tindel~\cite{tindelnourdin} studied asymptotic expansions for solutions of SDE driven by fractional Brownian motion using expansion over trees.

In Lyons' theory  to perform various
computations (e.g. Taylor expansions) the geometric
condition is (implicitly) used to  ensure that products of iterated integrals can be
expanded in a sum of other iterated integrals. On the other hand
iterated integrals indexed by trees already form a closed algebra with
respect to point-wise product and path integration (see below for
details). Thus, by enriching the notion of rough path we are  able to perform
computations as in the case of geometric rough paths and build a
complete theory for non-geometric rough integrals. 
Moreover we hope that such a bridge can inspire novel integration
methods for stochastic differential equations in the line of~\cite{MR2198539}.

\bigskip
The plan of the note is the following. In Sect.~\ref{sec:trees} we
introduce the concept of (labeled) rooted tree, the associated
(D\"urr-Connes-Kreimer) Hopf algebra and fix the relative notations. In
Sect.~\ref{incr} we summarize the theory of finite increments
described in~\cite{MR2091358} which can be used as the base for
building  rough paths theory. In Sect.~\ref{sec:itin} we
introduce iterated integrals indexed by labeled rooted trees and
prove the basic multiplicative property which is a generalization of
Chen's multiplicative property for usual iterated integrals. Next,
in Sect.~\ref{sec:trees-diff} we explain how sums over iterated
integrals indexed by rooted trees encode  the
solutions of driven differential equations. At this point we are ready
to generalize rough paths and introduce the notion of \emph{branched
  rough path} (in Sect.~\ref{sec:brp}), prove a generalized
extension theorem and construct  the branched rough path
associated to an \emph{ almost} branched rough path (following the
development of the standard theory, see
e.g.~\cite{MR1654527}). In Sect.~\ref{sec:controlled}, we
introduce path controlled by a branched rough path and show how to solve
differential equations driven by a branched rough path.
Finally in Sect~\ref{sec:inf-dim} we discuss another motivation to
consider tree-labeled series: rough paths adapted to the solution of
infinite-dimensional equations (deterministic or stochastic).

\section{Trees}
\label{sec:trees}
Given a finite set $\cl$, define a $\cl$-labeled rooted tree as a
finite graph with a special vertex called \emph{root} such
that there is a unique path from the root to any other vertex of the
tree. Moreover to each vertex there is associated an element of $\cl$. Here some examples of rooted trees labeled by $\cl=\{1,2,3\}$:
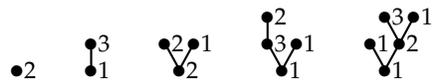
$$
\tnodel 2 \qquad 
\pstree{\tnodel 1}{\tnodel 3}
\qquad
\pstree{\tnodel 2}{\tnodel 2 \tnodel 1}
\qquad
\pstree{\tnodel 1}{\pstree{\tnodel 3}{\tnodel 2} \tnodel 1}
\qquad
\pstree{\tnodel 1}{\tnodel 1 \pstree{\tnodel 2}{\tnodel 3 \tnodel 1}} 
$$
We draw the root  at the bottom with the tree growing upwards.
Note that in a rooted tree the order of the branches at any vertex is ignored so
the following two are representations of the same (unlabeled) tree:
$$
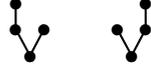
\pstree{\tnode}{\pstree{\tnode}{\tnode} \tnode}
\qquad
\pstree{\tnode}{\tnode \pstree{\tnode}{\tnode} }
$$
Given $k$ $\cl$-decorated rooted trees $\tau_1,\cdots,\tau_k$ and a label $a \in\cl$ we define
$\tau = [\tau_1,\cdots,\tau_k]_a$ as the tree obtained by attaching the $k$
roots of $\tau_1,\cdots,\tau_k$ to a new vertex with label $a$ which will be the root
of $\tau$. Any decorated rooted tree can be constructed using the
simple decorated tree $\bullet_a$ ($a \in\LL$)
and the operation $[\cdots]$, e.g.
$$
[\bullet] = \pstree{\tnode}{ \tnode}
\qquad
[\bullet,[\bullet]] = \pstree{\tnode}{\pstree{\tnode}{\tnode} \tnode},
\qquad \text{etc\dots}
$$

Denote $\TT_\cl$ the set of all $\cl$
decorated rooted trees and let $\TT$ the set of rooted trees without
decoration (i.e. for which the set of labels $\cl$ is made of a single
element). There is a canonical map $\TT_\cl \to \TT$ which simply
forget all the labels and every function on $\TT$ can be extended,
using this map to a function on $\TT_\cl$ for any set of labels
$\cl$. 
Let $|\cdot| : \TT \to \RR$ the map which counts the
number of vertices of the (undecorated) tree and which can be defined
recursively as
$$
|\bullet | = 1, \qquad |[\tau_1,\dots,\tau_k]| = 1+|\tau_1|+\cdots+|\tau_k|
$$
moreover we define the \emph{tree factorial} $\gamma: \TT \to \RR$ as
$$
\gamma(\bullet) = 1, \qquad \gamma([\tau_1,\dots,\tau_k]) =
|[\tau_1,\dots,\tau_k]| \gamma(\tau_1) \cdots \gamma(\tau_k)
$$

Last we define the \emph{symmetry factor} $\sigma : \TT_\cl \to \RR$ with the
recursive formula $\sigma(\tau) = 1$ for $|\tau| = 1$ and
\begin{equation}
  \label{eq:sigma-prop}
\sigma([\tau^1 \cdots \tau^k]_a) = \frac{k!}{\delta(\tau^1,\dots,\tau^k)} \sigma(\tau^1) \cdots \sigma(\tau^k)  
\end{equation}
where $\delta(\tau^1,\cdots,\tau^{k})$ counts the number of different ordered $k$-uples $(\tau^1,\cdots,\tau^{k})$ which corresponds to the same (unordered) collection  $\{\tau^1,\cdots,\tau^{k}\}$ of subtrees.
The factor $k!/\delta(\tau^1,\dots,\tau^k)$ counts the order of the
subgroup of permutations of $k$ elements which does not change the
ordered $k$-uple $(\tau^1,\cdots,\tau^{k})$. Then $\sigma(\tau)$ is
 is the order of the subgroup of permutations on the vertex of
the tree $\tau$ which do not change the tree (taking into account also
the labels). Another equivalent recursive definition for $\sigma$ is
$$
\sigma([(\tau^1)^{n_1}\cdots(\tau^k)^{n_k}]_a) = n_1!\cdots n_k!
\sigma(\tau^1)^{n_1}\cdots \sigma(\tau_k)^{n_k}
$$
where $\tau^1,\dots,\tau^k$ are distinct subtrees and $n_1,\dots,n_k$
the respective multiplicities.

Define the algebra $\cA \TT_\cl$ as the commutative polynomial algebra
generated by $\{1\}\cup \TT_\cl$ over $\RR$, i.e. elements of $\cA \TT_\cl$ are
finite linear combination with coefficients in $\RR$ of formal monomials in the form $\tau_1
\tau_2 \cdots \tau_n$ with $\tau_1,\dots,\tau_n \in \TT_\cl$ or of the
unit $1\in \cA\TT_\cl$. The set of all tree monomials is the set of
\emph{forests} $\FF_\cl$ including the empty forest $1 \in \FF_\LL$. The algebra $\cA \TT_\cl$
 is endowed with a graduation $g$ given  by $g(\tau_1
\cdots \tau_n) = |\tau_1|+\cdots+|\tau_n|$ and $g(1) = 0$. This graduation induces a
corresponding filtration of $\cA \TT_\cl$ in finite dimensional linear
subspaces $\cA_n \TT_\cl$ generated by the set $\FF_\cl^n$ of forests of degree $\le n$.

Any map $f : \TT_\cl \to
A$ where $A$ is some commutative algebra, can be extended in a unique way to
a homomorphism $f : \cA
\TT_\cl \to A$ by setting:
$
f(\tau_1 \cdots \tau_n) = f(\tau_1) f(\tau_2) \cdots  f(\tau_n)
$.

On the algebra $\cA \TT_\cl$ we can define a counit $\eps: \cA
\TT_\cl \to \RR$ as an algebra homomorphism such that $\eps(1) = 1$
and $\eps(\tau) = 0$ otherwise and a \emph{coproduct} $\Delta :
\cA\TT_\cl \to \cA\TT_\cl \otimes \cA\TT_\cl$ in the
following  way: $\Delta$ is an algebra homomorphism,
i.e. $\Delta(1) = 1 \otimes 1$, $\Delta(\tau_1 \cdots \tau_n) =
\Delta(\tau_1) \cdots \Delta(\tau_n)$ and acts linearly on linear
combinations of forests and on each tree it acts recursively as
\begin{equation}
  \label{eq:delta-def}
\Delta(\tau) = 1\otimes \tau  + \sum_{a \in \cl} (B^a_+ \otimes \text{id} )[\Delta(B^a_-(\tau))] 
\end{equation}
where $B_+^a(1) = \troot_a$ and $B_+^a(\tau_1 \cdots \tau_n) = [\tau_1 \cdots \tau_n]_a$ and
$B_-^a$ is the inverse of $B_+^a$ or is equal to zero if the tree
root does not have label $a$, i.e. 
\begin{equation*}
B_-^a(B^b_+(\tau_1 \cdots \tau_n))= \begin{cases}
 \tau_1 \cdots \tau_n & \text{if $a=b$}\\
0 & \text{otherwise}
\end{cases}
\end{equation*}
The coproduct $\Delta$ has an explicit description in terms of cuts
which is useful in some proofs. A \emph{cut} of a tree $\tau$ is a
subset of its edges which is selected to be removed. A cut is
\emph{admissible} if going from the root to any leaf of the tree we
meet at most one cut. Given a tree $\tau\in\TT_\LL$ and an admissible
cut $c$, we denote with $R_c(\tau)\in \TT_\LL$ the  tree obtained after the cut (that is the subgraph
containing the root) while the set
of subtrees detached from the ``trunk'' by the cut is denoted by
$P_c(\tau) \in \FF_\LL$. With this notation the action of the
coproduct on trees $\tau \in \TT_\LL$ can be described by the formula 
\begin{equation}
  \label{eq:copr-cuts}
  \Delta(\tau) = 1 \otimes \tau + \tau \otimes 1 + \sum_c P_c(\tau)\otimes R_c(\tau)
\end{equation}
where the sum is performed over all the admissible cuts $c$ of $\tau$.

Endowed with $\eps$ and $\Delta$ the algebra $\cA \TT_\cl$ become a
bialgebra, there exists also an antipode $S$ which complete the
definition of the Hopf algebra structure on $\TT_\cl$ as described by
Connes-Kreimer~\cite{MR1660199} (in the unlabeled case). 

Note that our definition of the coproduct differ from the one commonly
present in the literature by the exchange of the order of the factors
in the tensor product in order to be consistent with other notations
present in the paper. 

There exists various notations for the coproduct $\Delta$ we will often use Sweedler's notation $\Delta \tau = \sum \tau_{(1)} \otimes \tau_{(2)}$ but we also introduce a counting function $c : \TT_\cl \times \TT_\cl \times \FF_\cl \to \NN$ such that
$$
\Delta \tau = \sum_{\rho \in \TT_\cl, \sigma \in \FF_\cl} c(\tau,\rho, \sigma) \rho \otimes \sigma.
$$

In the following we will use letters $\tau,\rho,\sigma,\dots$ to
denote trees in $\TT_\cl$ or forests  in $\FF_\LL$, the degree $g(\tau)$ of a forest $\tau \in \FF_\cl$ will also be written as $|\tau|$. Roman letters $a,b,c,\dots \in \LL$ will denote vector  indexes (i.e. labels) while $\ol{a},\ol{b},\dots$ will denote multi-indexes with values in $\LL$: $\overline{a}=(a_1,\dots,a_n) \in \LL^n$ with $|\overline a|=n$ the size of this multi-index.

\section{Increments}
\label{incr}
Given $T > 0$, a vector space $V$ and an integer $k \ge 1$, we denote by
$\CC_k(V)$ the set of functions $g : [0,T]^{k} \to V$
 such
that $g_{t_1 \cdots t_{k}} = 0$
whenever $t_i = t_{i+1}$ for some $0 \le i \le k-1$.
Such a function will be called a 
\emph{$k$-increment}, and we will
set $\CC_*(V)=\cup_{k\ge 1}\CC_k(V)$. We write $\CC_k =
\CC_k(\RR)$. There is a cochain complex $(\CC_*(V),\delta)$ where the coboundary $\delta$, satisfying $\delta^2 = 0$, is defined as follows on $\CC_k(V)$:
\begin{equation}
  \label{eq:coboundary}
\delta : \CC_k(V) \to \CC_{k+1}(V) \qquad 
(\delta g)_{t_1 \cdots t_{k+1}} = \sum_{i=1}^{k+1} (-1)^i g_{t_1
  \cdots \hat t_i \cdots t_{k+1}} , 
\end{equation}
here $\hat t_i$ means that this particular argument is omitted. 
 We will denote $\cZ\CC_k(V) = \CC_k(V) \cap \text{Ker}\delta$ 
and $\cB \CC_k(V) =
\CC_k(V) \cap \text{Im}\delta$, respectively the spaces of
\emph{$k$-cocycles} and of \emph{$k$-coboundaries}.

Some simple examples of actions of $\der$,
which will be the ones we will really use throughout the paper,
 are obtained by letting 
$g\in\CC_1(V)$ and $h\in\CC_2(V)$. Then, for any $t,u,s\in\ott$, we have
$
  (\der g)_{ts} = g_t - g_s $, and $
(\der h)_{tus} = h_{ts}-h_{tu}-h_{us}
$.
Furthermore, it is readily checked~\cite{MR2091358} that
the complex $(\CC_*(V),\delta)$ is \emph{acyclic}, i.e. 
$\cZ \CC_{k+1}(V) = \cB \CC_{k}(V)$ for any $k\ge 1$, or otherwise stated, 
the sequence
\begin{equation}
\label{eq:exact_sequence}
 0 \rightarrow \RR \rightarrow \CC_1(V)
 \stackrel{\der}{\longrightarrow} \CC_2(V) \stackrel{\der}{\longrightarrow} \CC_3(V) \stackrel{\der}{\longrightarrow} \CC_4(V) \rightarrow \cdots 
\end{equation}
is exact. 
This implies in particular that if $\der h = 0$ for some $h \in \CC_{2}(V)$ then there exists $f \in \CC_{1}(V)$ such that $\der f = h$.
 Thus we get a heuristic 
interpretation of the coboundary $\der h$:  it measures how much a
given 2-increment $h$  is far from being an {\it exact} increment of  a
function (i.e. a finite difference). 

\bigskip
When  $V=\R$ the complex $(\CC_*,\delta)$ is an (associative, non-commutative)
graded algebra once endowed with the following (exterior) product:
for  $g\in\CC_n $ and $h\in\CC_m $ let  $gh \in \CC_{n+m-1} $
the element defined by
\begin{equation}\label{cvpdt}
(gh)_{t_1,\dots,t_{m+n-1}}=
g_{t_1,\dots,t_{n}} h_{t_{n},\dots,t_{m+n-1}},
\quad
t_1,\dots,t_{m+n-1}\in\ott.
\end{equation}
In this context,
the coboundary $\delta$ act as a graded derivation with respect to the
algebra structure. In particular we have the following useful properties.
\begin{enumerate}
\item
Let $g,h$ be two elements of $\CC_1 $. Then
\begin{equation}\label{difrulu}
\der (gh) = \der g\,  h + g\, \der h.
\end{equation}
\item
Let $g \in \CC_1 $ and  $h\in \CC_2 $. Then
$$
\der (gh) = \der g\, h + g \,\der h, \qquad
\der (hg) = \der h\, g  - h \,\der g.
$$
\end{enumerate}
\bigskip

The iterated integrals of smooth functions on $\ott$ are 
particular cases of elements of $\CC$ which will be of interest for
us. Consider $f\in\CC_1^\infty $, where $\CC_1^\infty $ is the set of
smooth functions from $\ott$ to $\R$. For each $h \in \CC_2$ the integral $\int_s^t df_u \,
h_{us}$, which will be denoted by  
$\cj(df \,  h)$, can be considered as an element of
$\CC_{2}$. That is, for $s,t\in\ott$, we set 
$$
\cj_{ts}(df \,  h)
= \int_s^t  dg_u h_{us}.
$$

The basic relation between integration and the coboundary $\der$ is given by the next lemma.

\begin{lemma}
\label{lemma:split}
Let $h \in \CC_2$ such that $\der h = \sum_i h^{1,i} h^{2,i}$ (finite sum) for $h^{1,i},h^{2,i} \in \CC_2$ and let $x \in \CC^\infty_1$. Then
\begin{equation}
  \label{eq:action-der}
\der \cj(dx\, h) = \cj(dx) h + \sum_i \cj(dx\, h^{(1,i)}) h^{(2,i)}  
\end{equation}
\end{lemma}
\begin{proof}
  \begin{equation*}
    \begin{split}
  \der \cj(dx\, h)_{tus} & = \int_s^t h_{vs} dx_v - \int_s^u h_{vs} dx_v - \int_{u}^t h_{vu}dx_v
\\ & = \int_u^t (h_{vs}-h_{vu}) dx_v
 = \int_u^t \der h _{vus} dx_v + \int_{u}^t h_{us} dx_v
\\ & = \sum_i \int_u^t h^{(1,i)}_{vu}   dx_v\, h^{(2,i)}_{us} + \cj_{tu}(dx)\, h_{us}
    \end{split}
  \end{equation*}
\end{proof}

Then given a vector $\{x^i \}_{i=1,\dots,d}$ of elements of $\CC_1^{\infty}$ introduce iterated integrals
 recursively as
$$
\cj(dx^{i_1} dx^{i_2}  \cdots dx^{i_n}) = \cj[dx^{i_1} \cj(dx^{i_2} \cdots dx^{i_n})].
$$ 
where $i_1,\dots,i_n \in \{1,\dots,d\}$.
Then by using Lemma~\ref{lemma:split} we recover Chen's multiplicative property (in disguise)
\begin{equation}
  \label{eq:chen}
\der \cj(dx^{i_1} \cdots dx^{i_n}) = \sum_{k=1}^{n-1} \cj(dx^{i_{1}} \cdots dx^{i_k}) \cj(dx^{i_{k+1}}\cdots dx^{i_n}), \qquad (i_1,\dots,i_n) \in \{1,\dots,d\}^n.
\end{equation}

\section{Rooted trees and iterated integrals}
\label{sec:itin}
Fix a family $x= \{x^a\}_{a=1,\dots,d}$ of smooth elements in $\CC_1$
and let $\cl =\{1,2,\dots,d\} $ the set of indexes. 

By iterating integrations along the elements of $x$ we can build a map $X :
\TT_\cl  \to C([0,T]^2; \RR)$ defined as follows
\begin{equation}
\label{eq:it-integrals}
X_{ts}^{\troot_a} = \int_s^t dx_u^a,
\qquad
X_{ts}^{[\tau^1 \cdots \tau^k]_a} = \int_s^t 
\prod_{i=1}^k X_{us}^{\tau^i}
dx_u^a. 
\end{equation}
On the vector space $\cac_2$ we introduce
the associative and commutative inner product $\circ$ as $(a \circ
b)_{ts} = a_{ts} b_{ts}$ for $a,b \in \cac_2$. With this product
$\cac_2$ becomes an algebra and as explained before we can extend the
map $X : \TT_\cl \to \cac_2$ to a map on $\cA
\TT_\cl$ by linearity and by letting $X^{\tau_1 \cdots \tau_n}_{ts} =
X^{\tau_1}_{ts} X^{\tau_2}_{ts} \cdots X^{\tau_n}_{ts}$ for the value
of $X$ on the forest $\tau_1\cdots\tau_n$. Using this product we can
write 
$X^{[\tau_1\cdots\tau_n]_a} = \int X^{\tau_1 \cdots \tau_n} dx^a$.

Let $\CC^+_2 = \CC_2 \oplus e$ the unital  algebra obtained by adding
to the algebra $\CC_2$ the unit $e$ such that $e_{ts} = 1$ for any
$t,s\in[0,T]$. 

The
product $\circ$ has the following relation with $\der$:
\begin{equation}
 \label{eq:der-leib-2}
  \der (a \circ b) = \der a \circ \der b + (e a + a e)\circ \der b +
(eb+be) \circ \der a + ab + ba
\end{equation}
 where
$\circ$ is defined on $\cac_3$ in the natural way: $(g \circ
h)_{tus} = g_{tus} h_{tus}$ for every $g,h \in \cac_3$.

If on the algebra $(\cac_2,\circ)$ we consider the  exterior product $\cac_2 \otimes \cac_2
\to \cac_3$ then we can extend the homomorphism $X$ also to the tensor
product $\cA \TT_\cl \otimes \cA \TT_\cl$ by $X^{\sigma \otimes \rho}
= X^{\sigma}X^{\rho}$ for every $\sigma,\rho \in \cA \TT_\cl$.

Denote with $I^a : \CC_2 \to \CC_2$ the integration map given
by $I^a(h) = \cj( dx^a h)$ then for all elements $\sigma \in
\cA\TT_\cl$ we have  $I^a X^\sigma = X^{B_+^a \sigma}$: the map
$B_+^a$ represent integration  on the sub-algebra $\cA_X \subset
\CC^+_2$ generated by $\{X^\tau\}_{\tau \in \TT_\LL}$.
This sub-algebra contains the polynomial algebra  generated by the set $\{\der x^a\}_{a\in\LL}$:
\begin{equation}
  \label{eq:pol-in-trees}
X^{\troot_{a_1} \cdots \troot_{a_n}} = X^{\troot_{a_1}} \circ \cdots \circ X^{\troot_{a_n}}
= \der x^{a_1} \circ \cdots \circ  \der x^{a_n}.  
\end{equation}
It contains also the usual iterated integrals of $x$:
\begin{equation}
  \label{eq:chen-in-trees}
\cj(dx^{a_1} \cdots dx^{a_n}) = I^{a_1} I^{a_2} \cdots I^{a_{n-1}}( \der x^{a_n}) = X^{B_+^{a_1} B_+^{a_2}\cdots B_+^{a_{n-1}} \troot_{a_n}}
=
 X^{[\cdots[\troot_{a_n}]_{a_{n-1}}\cdots]_{a_1}}.  
\end{equation}
To future use let us denote with $\TT^{\text{Chen}}_\cl$ the subset of $\TT_\cl$ made of ``linear'' labeled trees of the form $[\cdots[\troot_{a_n}]_{a_{n-1}}\cdots]_{a_1}$.

What is remarkable is the relation between the
coalgebra structure of the trees and  the algebraic properties of the
iterated integrals $X$ with respect to the coboundary $\delta$ as illustrated in the next theorem.

\begin{theorem}[Tree multiplicative property]
\label{th:equiv}
The map $X$ satisfy
the following algebraic relation:
\begin{equation}
  \label{eq:alg-relations}
\der X^{\sigma} = X^{\Delta'(\sigma)} , \qquad \sigma \in \cA\TT_\cl
\end{equation}
where $\Delta'$ is the reduced
coproduct
$
\Delta'(\tau) = \Delta(\tau) - 1 \otimes \tau - \tau \otimes 1.
$  
\end{theorem}

\begin{proof}
We will proceed by induction on the degree $g$  of the forests in $\cA
\TT_\cl$ defined above. It is clear that the
relation~(\ref{eq:alg-relations}) holds for the simple tree $\troot_a$
with degree $g=1$. Assume that eq.~(\ref{eq:alg-relations}) holds for
every monomial with degree less than $n$ and let us prove it for
monomials of degree $n$.

We need the following two properties of the reduced coproduct: first, its
recursive definition can be rewritten as
\begin{equation}
  \label{eq:delta-prime-rec}
 \Delta'(\tau) = \sum_{a\in\LL} \troot_a \otimes B^a_-(\tau)  + \sum_a (B_+^a
 \otimes \text{id}) [\Delta'(B_-^a(\tau))]    
\end{equation}
which follows directly from~(\ref{eq:delta-def}), next a formula for
the action of $\Delta'$ on products of monomials:
\begin{equation}
\label{eq:leib-delta-prime}
\Delta'(\rho \sigma) = \Delta' \sigma \Delta' \rho 
+ (1\otimes \sigma + \sigma \otimes 1) \Delta' \rho 
+ (1\otimes \rho + \rho \otimes 1) \Delta' \sigma 
+ \rho \otimes \sigma + \sigma \otimes \rho.  
\end{equation}
for  $\rho,\sigma$ monomials on trees. Assume $g(\rho\sigma) = n$ and let us compute $\der X^{\rho
  \sigma}$ using eq.~(\ref{eq:der-leib-2}):
\begin{equation*}
  \begin{split}
\der X^{\rho \sigma} & = \der (X^\rho \circ X^\sigma) 
\\ & =   
\der X^\rho \circ (X^\sigma e + e X^\sigma) + \der X^\sigma \circ (X^\rho e + e
X^\rho) + \der X^\rho \circ \der X^\sigma + X^\rho X^\sigma + X^\sigma X^\rho    
  \end{split}
\end{equation*}
Since $g(\sigma) < n$ and $g(\rho) < n$ we obtain
\begin{equation*}
  \begin{split}
\der X^{\rho \sigma}
& =   
X^{\Delta' \rho} \circ (X^\sigma e + e X^\sigma) +  X^{\Delta'\sigma} \circ (X^\rho e + e
X^\rho) +  X^{\Delta'\rho} \circ  X^{\Delta'\sigma} + X^\rho X^\sigma + X^\sigma X^\rho    
\\& =   
X^{\Delta' \rho} \circ X^{\sigma\otimes 1 + 1 \otimes \sigma} +
X^{\Delta'\sigma} \circ X^{\rho\otimes 1 + 1 \otimes \rho} +
X^{\Delta'\rho} \circ  X^{\Delta'\sigma} + X^{\rho \otimes \sigma} +
X^{\sigma \otimes \rho}
\\& =   
X^{\Delta' \rho( \sigma\otimes 1 + 1 \otimes \sigma)} +
X^{\Delta'\sigma (\rho\otimes 1 + 1 \otimes \rho)} +
X^{\Delta'\rho \Delta'\sigma} + X^{\rho \otimes \sigma} +
X^{\sigma \otimes \rho}
\\ & = X^{\Delta'(\rho \sigma) }
  \end{split}
\end{equation*}
according to eq.~(\ref{eq:leib-delta-prime}).
So we have proven eq.~(\ref{eq:alg-relations}) for nontrivial monomials
of $g$-degree $n$. It remains to prove the relation for monomials
given by a single tree of degree $n$. To do this we need the action of
$\der$ on iterated integrals which is given by Lemma~\ref{lemma:split} above.
Let us compute $\der X^\tau$ using formula~(\ref{eq:action-der})
with $\tau = [\tau_1 \cdots \tau_n]_a$:
\begin{equation*}
  \begin{split}
\der X^{[\tau_1 \cdots \tau_n]_a} & = \der \cj[ dx^a X^{\tau_1 \cdots \tau_n}]
=  \der x^a   X^{\tau_1 \cdots \tau_n} + \sum_i \cj[dx^a X^{\theta^1_i}]\, X^{\theta^2_i}
\\ & =    X^{\troot_a} X^{\tau_1 \cdots \tau_n} + \sum_i X^{[\theta^1_i]_a}\, X^{\theta^2_i}
  \end{split}
\end{equation*}
where $\der X^{\tau_1 \cdots \tau_n} = \sum_i X^{\theta^1_i}
X^{\theta^2_i}$ and $\theta^{1,2}$ satisfy
$
\Delta'(\tau_1 \cdots \tau_n) = \sum_i \theta^1_i \otimes \theta^2_i
$
since our induction assumptions imply that the monomial $\tau_1 \cdots \tau_n$,
eq.~(\ref{eq:alg-relations}) holds. 
Then
\begin{equation*}
  \begin{split}
\der X^{[\tau_1 \cdots \tau_n]_a}
 & =    X^{\troot_a \otimes  (\tau_1 \cdots \tau_n)} +
 X^{\sum_i [\theta^1_i]_a \otimes \theta^2_i}
=    X^{\troot_a \otimes  (\tau_1 \cdots \tau_n) +
 \sum_i [\theta^1_i]_a \otimes \theta^2_i}
\\ & =    X^{\troot_a \otimes  (\tau_1 \cdots \tau_n) +
  (B_+^a \otimes \text{id})(\sum_i \theta^1_i \otimes \theta^2_i)}
=    X^{\Delta'([\tau_1 \cdots \tau_n]_a)}
  \end{split}
\end{equation*}
where we used eq.~(\ref{eq:delta-prime-rec}). Then we proved eq.~(\ref{eq:alg-relations}).
\end{proof}

\begin{example}
\label{ex:example-1}
Let us give an example in one dimension ($d=1$) so trees are not decorated. The forests
of degree less or equal to three are:
\smalltrees
\begin{equation*}
\tsnode, \aabb, \tsnode \tsnode, \aaabbb, \tsnode \aabb, \tsnode \tsnode \tsnode, \aababb  
\end{equation*}
The reduced coproduct on these monomials acts as follows:
\begin{equation*}
\Delta' \aabb = \tsroot \otimes \tsroot,
\qquad \Delta' (\tsroot \tsroot) = 2 \tsroot \otimes \tsroot  
\end{equation*}
\begin{equation*}
\Delta'   \aaabbb = \aabb \otimes \tsroot + \tsroot \otimes \aabb 
\end{equation*}
\begin{equation*}
\Delta' ( \tsroot \aabb) = \tsroot \otimes \tsroot \tsroot + \tsroot \tsroot \otimes \tsroot
+ \aabb \otimes \tsroot + \tsroot \otimes \aabb
\end{equation*}
\begin{equation*}
\Delta'(\tsroot^3) = 3 \tsroot^2 \otimes \tsroot + 3 \tsroot \otimes \tsroot^2  
\end{equation*}
\begin{equation*}
 \Delta' \aababb  = \tsroot \otimes \tsroot \tsroot + 2 \aabb \otimes \tsroot
\end{equation*}

So we have
$$
\der X^{\aababb} = X^{\tsroot} X^{\tsroot \tsroot} + 2 X^{\aabb} X^{\tsroot}
$$
\end{example}

\begin{remark}
A particular case of the tree multiplicative property~(\ref{eq:alg-relations}) is given by Chen's multiplicative property~(\ref{eq:chen}) with the aid of the relation~(\ref{eq:chen-in-trees}).
\end{remark}

As a first elementary application of this result we derive a tree binomial formula.

\begin{lemma}[Tree Binomial]
\label{lemma:tree-binomial}
For every $\tau \in \TT$ and $a,b \ge 0$ we have
\begin{equation}
  \label{eq:tree-binomial}
  (a+b)^{|\tau|} = \sum_i \frac{\tau!}{\tau^{(1)}_i! \tau^{(2)}_i!} a^{|\tau^{(1)}_i|} b^{|\tau^{(2)}_i|}
\end{equation}
\end{lemma}
\begin{proof}
Consider the iterated integrals $T^\tau$ associated to the identity path $t : \RR \to \RR$
$$
T^{\troot}_{ts} = t-s, \qquad T^{[\tau_1 \cdots \tau_n]}_{ts} = \int_s^t T^{\tau_1}_{us}\cdots T^{\tau_n}_{us} du
$$  
By induction it is not difficult to prove that
$
T^{\tau}_{ts} = (t-s)^{|\tau|} (\tau!)^{-1}
$, so applying Thm.~\ref{th:equiv} to $T^\tau$ we get
\begin{equation*}
  \begin{split}
\frac{(t-s)^{|\tau|}}{\tau!} & = T^{\tau}_{ts}  = T^{\tau}_{us} + T^{\tau}_{tu} + \sum'_i T^{\tau^{(1)}_i}_{tu} T^{\tau^{(2)}_i}_{us}
= \sum_i T^{\tau^{(1)}_i}_{tu} T^{\tau^{(2)}_i}_{us}    
  \\ & = \sum_i \frac{1}{\tau^{(1)}_i! \tau^{(2)}_i!} (t-u)^{|\tau^{(1)}_i|} (u-s)^{|\tau^{(1)}_i|}
  \end{split}
\end{equation*}
Then setting $t-u = a$ and $u-s=b$ we get eq.~(\ref{eq:tree-binomial}).
\end{proof}

\subsection{Geometric paths}
The above homomorphism $X$ can be simplified using the fact that it is generated by a $C^1$ family $x$. Indeed Chen~\cite{MR0454968} proved that products of iterated integrals can be always expressed as linear combination of iterated integrals via the \emph{shuffle product}:
\begin{equation}
  \label{eq:shuffle}
\cj(dx^{a_1}\cdots dx^{a_n})\circ \cj(dx^{b_1}\cdots dx^{b_m}) = \sum_{\overline c \in \text{Sh}(\overline a, \overline b)} \cj(dx^{c_1}\cdots dx^{c_{n+m}})  
\end{equation}
where given two multi-indexes $\overline a = (a_1,\dots,a_n)$ and $\overline b = (b_1,\dots,b_n)$ their \emph{shuffles} $\text{Sh}(\overline a, \overline b)$ is the set of all the possible permutations of the $(n+m)$-uple $(a_1,\dots,a_n,b_1,\dots,b_m)$ which does not change the ordering of the two subsets $\overline a$, $\overline b$.

Using relation~(\ref{eq:shuffle}) we can reduce every $X^{\tau}$ for $\tau \in \TT_\cl$ to a linear combination of $\{ X^{\sigma} \}_{\sigma \in \TT^{\text{Chen}}_\cl}$.

\section{Series solutions of driven differential equations}
\label{sec:trees-diff}

Under appropriate conditions on the vectorfield $f: \RR^n \to \RR^n$
the solution $y$ of the differential equation $dy/dt = f(y)$ $y_0 =
\eta$ admit the series representation
\begin{equation}
  \label{eq:B-series}
y_t = \eta + \sum_{\tau \in \TT} \psi^f(\tau)(\eta)
\frac{t^{|\tau|}}{\sigma(\tau) \tau!}  
\end{equation}
which is called $B$-series (in honor of J.~Butcher,
see\cite{MR1993957,MR0403225,MR1227985}). The coefficients $\psi^f$ are called
\emph{elementary differentials} and are defined as
$$
\psi^f(\troot)(\xi) = f(\xi), \qquad \psi^f([\tau^1 \cdots \tau^k]) =
\sum_{\ol b \in \II\LL_1} f_{\ol b}(\eta)
\psi^f(\tau^1)(\xi)^{b_1}\cdots \psi^f(\tau^k)(\xi)^{b_k} 
$$ 
where we introduce multi-indexes $\overline b \in \II\LL_1 = \cup_{k=0}^\infty \LL_1^k$, $\LL_1 = \{1,\dots,n\}$, with the convention $\LL_1^0 = \emptyset$ and we set  $f_{\emptyset}(\xi) = f(\xi)$ and
$
f_{\ol b}(\xi) = \prod_{i=1}^{|\ol b|} \partial_{\xi_{b_i}} f(\xi)
$ for the derivatives of the vectorfield.
 
In this section we study the analogous series expansion for \emph{driven}
differential equation. 
Consider a $C^1$ path $x: [0,T] \to \RR^d$ and let $\{x^a\}_{a\in \LL}$ be its coordinates in a fixed basis. Fix a point $\eta \in \RR^n$ and let $f_a: \RR^n \to \RR^n$, $a=1,\dots,d$ be a collection of analytic vectorfields on $\RR^n$. Let $R$ be a common analiticy radius around $\eta$ for all coordinates.
\begin{theorem}
 The solution of the differential equation
$
 dy_t = \sum_{a\in\LL}  f_a(y_t) dx_t^a$, $y_0 = \eta
$  
admit locally the series representation
\begin{equation}
\label{eq:tree-rep}
\der y_{ts} = †\sum_{\tau \in \TT_\LL} \frac{1}{\sigma(\tau)} \phi^{f}(\tau)(y_s) X^\tau_{ts}, \qquad y_0 = \eta
\end{equation}
where the sum runs over all $\LL$-labeled rooted trees
$\tau\in\TT_\LL$
and where  we recursively define
functions $\phi^f : \TT_\LL \times \RR^n \to \RR^n$  such that
\begin{equation*}
\phi^f(\troot_a)(\xi) = f_a(\xi), \qquad
\phi^f([\tau^1 \cdots \tau^k]_a)(\xi) =  \sum_{\ol b  \in \II \LL_1: |\ol b| = k}
f_{a;b_1\dots b_k}(\xi) \prod_{i=1}^{k} [\phi^{f}(\tau^i)(\xi)]^{b_i}  .
\end{equation*}
\end{theorem}

\begin{proof}
Let us assume for the moment that the series~(\ref{eq:tree-rep}) converges absolutely.
We will verify that that eq.~(\ref{eq:tree-rep})  satisfy the integral
equation
\begin{equation}
\label{eq:integral-eq}
\der y_{ts}  = \sum_{a\in\LL} \int_s^t f_a(y_u) dx^a_u .  
\end{equation}

Consider the Taylor series for $f$ around $\xi\in \RR^n$:
$$
f_a(\xi')  = \sum_{\ol b\in \II\LL_1 } \frac{f_{a; \ol b}(\xi)}{|\ol b|!} \prod_{i=1}^{|\ol b|}
(\xi'-\xi)^{b_i} 
$$
where $\xi^k$ is the $k$-th coordinate of the vector $\xi\in \RR^n$. By the analyticity of the vectorfields $f_a$ this series converges as long as $|\xi-\xi'| \le R-|\xi'-\eta|$.

Compute the r.h.s. of eq.~(\ref{eq:integral-eq}) by plugging in
eq.~(\ref{eq:tree-rep}) and the Taylor expansion of $f$:
\begin{equation*}
  \begin{split}
 \sum_{a\in\LL} & \int_s^t f_a(y_u) dx_u^a  = \sum_{a\in\LL} \sum_{\ol b \in \II \LL_1} 
\frac{f_{a;\ol b}(y_s)}{|\ol b|!} \int_s^t \left( \prod_{i=1}^{|\ol b|}
\der y^{b_i}_{us}\right) dx_u^a   
\\ & = \sum_{a\in\LL} \sum_{\ol b \in \II \LL_1} 
\frac{f_{a;\ol b}(y_s)}{|\ol b|!} \int_s^t \prod_{i=1}^{|\ol b|}
\left[
\sum_{\tau \in \TT_\LL}  \frac{1}{\sigma(\tau)} [\phi^{f}(\tau)(y_s)]^{b_i} X^{\tau}_{us}
 \right] dx_u^a    
\\ & = \sum_{a\in\LL} \sum_{\ol b \in \II\LL_1} \frac{f_{a;\ol b}(y_s)}{|\ol b|!}  \sum_{\tau^1,\cdots,\tau^{|\ol b|}}
 \frac{1}{\sigma(\tau^1)\cdots \sigma(\tau^{|\ol b|})}\left(\prod_{i=1}^{|\ol b|}[ \phi^{f}(\tau^i)(y_s)]^{b_i}\right)
 \int_s^t \prod_{i=1}^{|\ol b|}
  X^{\tau^i}_{us}
 dx_u^a    
\\ & = \sum_{a\in\LL} \sum_{k=0}^\infty  \sum_{\tau^1,\cdots,\tau^{k}} \frac{1}{k!\sigma(\tau^1)\cdots \sigma(\tau^k)}\sum_{\ol b \in \II\LL_1 : |\ol b| = k} f_{a;\ol b}(y_s)  
 \left(\prod_{i=1}^{k}[\phi^{f}(\tau)(y_s)]^{b_i}\right)
 \int_s^t \prod_{i=1}^{k}
  X^{\tau^i}_{us}
 dx_u^a    
\\ & =\sum_{a\in\LL} \sum_{k=0}^{\infty} 
\sum_{\tau^1,\cdots,\tau^{k}}
\frac{1}{\sigma([\tau^1\cdots\tau^{k}]_a) \delta(\tau^1,\cdots,\tau^{k})}
\phi^f([\tau^1\cdots\tau^{k}]_a)(y_s) X_{ts}^{[\tau^1\cdots\tau^{k}]_a}
\\ & = \sum_{\tau \in \TT_\LL} \frac{1}{\sigma(\tau)} \phi^{f}(\tau)(y_s) X^\tau_{ts}
  \end{split}
\end{equation*}
which proves the claim . Note the multiplicity factor $\delta$ which disappears from the last line.

To prove the absolute convergence of the series we need bounds on $X^\tau$ and $\phi^f(\tau)$. For  $X^\tau$ we have:
$$
|X^{\tau}_{ts}| \le \frac{[A|t-s|]^{|\tau|}}{\tau!}
$$
where $A = \sup_{t\in[0,T]} |\dot x_t|$. This bound can be easily proven inductively on $\tau$.

Since $f_a$ are analytic functions, from Cauchy inequalities we obtain
$$
|f_{a,\ol b}(y_s)| \le \theta(\ol b) M  (R-r_s)^{-|\ol b|} \le |\ol b|! M (R-r_s)^{-|\ol b|} \le g^{(|\ol b|)}(r_s)
$$
see e.g.~\cite[pag. 47]{MR1227985}.
where $r_s = |y_s-\eta|$ and $M$ is a constant depending only on $\{f_a\}_{a\in \LL}$ and where we introduced the function $g(r) = M R (R-r)^{-1}$ and its derivatives $g^{(k)}(r) = MR k! (R-r)^{-k-1}$.
Define ``elementary differentials'' $\psi:\TT \times [0,R) \to \RR$ for $g$ as
$$
\psi(\troot)(r) = g(r), \qquad \psi([\tau_1 \cdots \tau_k])(r) = g^{(k)}(r) k! M (R-r)^{-k} 
$$
Then we have the bounds $|\phi^f(\tau)(y_s)| \le \psi(\tau)(r_s)$ for
any $\tau \in \TT_\LL$ and the series~(\ref{eq:tree-rep}) can be
bounded by
$$
 \sum_{\tau \in \TT_\LL} \frac{1}{\sigma(\tau)} \psi(\tau)(r_s) A^{|\tau|} \frac{|t-s|^{|\tau|}}{\tau!} 
$$
and by taking into account the multiplicity $d^{|\tau|}$ of labeled
trees corresponding to the same tree  $\tau$ we get 
$$
\sum_{\tau \in \TT} \frac{1}{\sigma(\tau)} \psi(\tau)(r_s) (dA)^{|\tau|} \frac{|t-s|^{|\tau|}}{\tau!} 
$$
This series is exactly the B-series~(\ref{eq:B-series}) for the solution $r_t$ of the differential equation 
\begin{equation}
\label{eq:g-eq}
\frac{dr_t}{dt} = dA g(r_t) = d A M R (R-r_t)^{-1}, \qquad r_0 = 0
\end{equation}
when written starting from $r_s$ at time $s< t$. 
Then 
$$
r_t = r_s+ \sum_{\tau \in \TT} \frac{1}{\sigma(\tau)} \psi(\tau)(r_s) (dA)^{|\tau|} \frac{|t-s|^{|\tau|}}{\tau!}
$$
as long as the solution $r_t$ exists and has a power series expansion in $t-s$. But the explicit solution of eq.~(\ref{eq:g-eq}) is given by
$r_t = R(1-\sqrt{1-t/t_*})$ with $t_* = R/(2dAM)$ and has power series expansion for any $t < t_*$.
So the original series is summable at least for any $t,s \in [0,t_*)$.  
\end{proof}

In the rest of this section we will denote $y^\tau_s = \phi^f(\tau)(y_s)/\sigma(\tau)$ so that
$
\der y_{ts} = \sum_{\tau \in \TT_\LL} X^\tau_{ts} y^\tau_s
$
moreover we will use the convention $X^{\emptyset}_{ts} = 1$ and $y^\emptyset_s = y_s$ to write
$$
y_t = \sum_{\tau \in \TT_\LL \cup \{ \emptyset \}} X^{\tau}_{ts} y^{\tau}_s
$$
The recursion for $y^\tau$ reads
\begin{equation}
  \label{eq:ytau-rec}
y^{\troot_a}_s = f_a(y_s), \qquad y^{[\tau^1 \cdots \tau^k]_a}_s = \frac{\sigma(\tau^1)\cdots \sigma(\tau^k)}{\sigma(\tau)} \sum_{\ol b: |\ol b|=k} f_{a, \ol b}(y_s) y_s^{\tau_1,b_1} \cdots y_s^{\tau_k,b_k}  
\end{equation}

We have the following theorem which show that each of the paths $y^\tau$ can be expanded in series w.r.t. to $X$ with coefficients which depends on the combinatorics of the reduced coproduct:

\begin{theorem}
\label{th:der-eqns-1}
For any $\tau \in \TT_\LL\cup \{ \emptyset \}$ we have
\begin{equation}
  \label{eq:ytau-der}
\der y^\tau_{ts} = \sum_{\sigma \in \TT_\LL, \rho \in \FF_\LL} c'(\sigma,\tau,\rho) X^{\rho}_{ts} y^{\sigma}_s  
\end{equation}
where $c'$ is the counting function for the reduced coproduct:
$
\Delta' \sigma = \sum_{\tau,\rho} c'(\sigma,\tau,\rho) \tau \otimes \rho
$.
\end{theorem}

\begin{proof}
The proof is by induction on $\tau$. The case $\tau = \troot_a$ requires only Taylor expansion:
\begin{equation}
  \label{eq:ytau-step-0}
  \begin{split}
\der y^{\troot_a}_{ts} & = \der f_a(y)_{ts} = \sum_{\ol b} \frac{f_{a; \ol b}(y)}{|\ol b|!} \left( \der y_{ts}\right)^{\ol b} 
\\  & = \sum_{k \ge 1} \sum_{\tau^1,\dots,\tau^{k}} \sum_{\ol b: |\ol b|=k} \frac{f_{a; \ol b}(y)}{k!}  y^{\tau^1,b_1}_s \cdots y^{\tau^k,b_k}_s X^{\tau^1 \cdots \tau^k}_{ts}
\\  & = \sum_{k \ge 1} \sum_{\tau^1,\dots,\tau^{k}} \frac{\sigma([\tau^1 \cdots \tau^k]_a)}{k! \sigma(\tau^1)\cdots\sigma(\tau^k)} y^{[\tau^1 \cdots \tau^k]_a}_s X^{\tau^1 \cdots \tau^k}_{ts}
\\  & = \sum_{k \ge 1} \sum_{\tau^1,\dots,\tau^{k}} \frac{1}{\delta(\tau^1,\dots,\tau^k)} y^{[\tau^1 \cdots \tau^k]_a}_s X^{\tau^1 \cdots \tau^k}_{ts}
\\  & = \sum_{\tau}  c'(\tau,\troot_a,\rho) y^{\tau}_s X^{\rho}_{ts}
  \end{split}
\end{equation}
since $c'(\tau,\troot_a,\rho) $ is different from zero, and take value one, iff $\tau = [\rho]_a$.

Now, assume eq.~(\ref{eq:ytau-der}) holds for all $\tau \in \TT^n_{\LL}$ and let us prove that it holds for trees $\tau$ with $|\tau|=n+1$. So take $\tau= [\tau^1 \dots \tau^k]_a$ with $|\tau|=n+1$, then $|\tau^i|\le n$ for any $i=1,\dots,k$. To compute the action of the map $\der$ on $y^\tau$ we use the recursive relation~(\ref{eq:ytau-rec}):
\begin{equation}
  \label{eq:ytau-step-2-der}
  \begin{split}
 \der y^{[\tau^1 \cdots \tau^k]_a}_{ts} & = \frac{\sigma(\tau^1)\cdots \sigma(\tau^k)}{\sigma(\tau)} \sum_{\ol b: |\ol b|=k} \der[ f_{a, \ol b}(y) y^{\tau^1,b_1} \cdots y^{\tau^k,b_k}]_{ts}  
   \end{split}
 \end{equation}
and the Leibniz formula
$$
\der (g^1 \cdots g^k)_{ts} = (g^1_s + \der g^1_{ts}) \cdots (g^1_s + \der g^1_{ts}) - g^1_s \cdots g^k_s = \sum_{G} G^1_{ts} \cdots G^k_{ts}
$$
where the sum is over all possible choices of $G$-s such that $G^i_{ts}= g^i_s$ or $G^i_{ts} = \der g^i_{ts}$ excluding the case where all the $G$-s are $g$ (that is, there should be at least one factor of the form $\der g^i$). By Taylor expansion
$$
\der f_{a, \ol b}(y)_{ts} = \sum_{m\ge 1} \sum_{\ol c: |\ol c|=m} \frac{f_{a,\ol b \ol c}(y)_s}{m!} \sum_{\eta^1,\dots,\eta^m} y^{\eta^1,c_1}_s\cdots y^{\eta^m,c_m}_s X_{ts}^{\eta^1\cdots \eta^m} 
$$
while using the induction hypothesis we have
$$
\der y^{\tau^i} = \sum_{\rho^i,\zeta^i} c(\zeta^i,\tau^i,\rho^i) X^{\rho^i} y^{\zeta^i} = \sum_{\zeta^i}  X^{\zeta^i_{(2)}} y^{\zeta^i} \delta_{\tau^i,\zeta^i_{(1)}}
$$
where there is an implicit sum over the terms $\zeta^i_{(1)},\zeta^i_{(2)}$ in the reduced coproduct of $\zeta^i$ and where $\delta_{\tau^i,\zeta^i_{(1)}}$ denotes the Kronecker delta function.
Then we rewrite eq.~(\ref{eq:ytau-step-2-der}) as
\begin{equation}
  \label{eq:ytau-step-2-2}
  \begin{split}
 \der y^{[\tau^1 \cdots \tau^k]_a}_{ts} & = \frac{\sigma(\tau^1)\cdots \sigma(\tau^k)}{\sigma(\tau)} 
\\ & \times \sum_{m\ge 0}  \frac{1}{m!} \sum_{\zeta^1,\dots,\zeta^k} \sum_{\eta^1,\dots,\eta^m}  
\sum_{\ol c: |\ol c|=m+k} f_{a,\ol c}(y_s)  y^{\eta^1,c_1}_s\cdots y^{\eta^m,c_m}_s  y^{\zeta^1,c_{m+1}}_s\cdots y^{\zeta^k,c_{m+k}}_s
\\ & \times X_{ts}^{\eta^1\cdots \eta^m\cdots \zeta^1_{(2)}\cdots \zeta^k_{(2)}} \delta_{\tau^1,\zeta^1_{(1)}} \cdots \delta_{\tau^k,\zeta^k_{(1)}}
   \end{split}
 \end{equation}
The summation in this formula has to be understood as follows:
the sum over $\zeta^i$ is performed on all trees which contains $\tau^i$ in the sense that $c'(\zeta^i,\tau^i,\rho^i)$ is different form zero for some $\rho^i$ and on the tree $\zeta^i=\tau^i$ in which case we understand that $\zeta^i_{(1)}=\tau^i$ and $\zeta^i_{(2)} = \emptyset$ (the empty forest). Note that this case in not contained in the reduced coproduct but is generated by the Leibniz's formula. Moreover we implicitly exclude from the summation above the case when $m=0$ and all the $\zeta^i$ are equal to the corresponding $\tau^i$. Then with this proviso we can simplify the above formula as
\begin{equation}
  \label{eq:ytau-step-2-3}
  \begin{split}
 \der y^{[\tau^1 \cdots \tau^k]_a}_{ts} & = \frac{\sigma(\tau^1)\cdots \sigma(\tau^k)}{\sigma(\tau)} 
\\ & \times \sum_{m\ge 0} \frac{1}{m!}\sum_{\zeta^1,\dots,\zeta^k} \sum_{\eta^1,\dots,\eta^m}  
\frac{\sigma(\zeta)}{\sigma(\zeta^1)\cdots\sigma(\zeta^k) \sigma(\eta^1)\cdots\sigma(\eta^m)}
 X^{\eta^1\cdots \eta^m\cdots \zeta^1_{(2)}\cdots \zeta^k_{(2)}} y^{\zeta} \delta_{\tau^1,\zeta^1_{(1)}} \cdots \delta_{\tau^k,\zeta^k_{(1)}}
   \end{split}
 \end{equation}
where $\zeta = [\zeta^1\cdots\zeta^k \eta^1\cdots \eta^k]$. Now,
recalling eq.~(\ref{eq:sigma-prop}), write
\begin{equation}
  \label{eq:ytau-step-2-4}
  \begin{split}
 \der y^{[\tau^1 \cdots \tau^k]_a}_{ts} & =  \sum_{m\ge 0} \sum_{\zeta^1,\dots,\zeta^k} \sum_{\eta^1,\dots,\eta^m}
\frac{(k+m)!}{k!m!}
\frac{\delta(\tau^1,\dots,\tau^k)}{\delta(\zeta^1,\dots,\zeta^k,\eta^1,\dots,\eta^k)}
 X^{\eta^1\cdots \eta^m\cdots \zeta^1_{(2)}\cdots
  \zeta^k_{(2)}} y^{\zeta} \delta_{\tau^1,\zeta^1_{(1)}} \cdots
\delta_{\tau^k,\zeta^k_{(1)}} .
   \end{split}
 \end{equation}
Introduce a new function $\tilde c: \TT_\LL \times \TT_\LL \times \FF_\LL \to \mathbb{N}$ such that
$$
\tilde c(\kappa_1,\kappa_2,\kappa_3) = \begin{cases}
c'(\kappa_1,\kappa_2,\kappa_3) & \text{for $\kappa_3 \neq \emptyset$}\\
\delta_{\kappa_1,\kappa_2} & \text{for $\kappa_3 = \emptyset$}
\end{cases}
$$
which counts the number of ways to cut away a forest $\kappa_3$ from the tree $\kappa_1$ leaving the tree $\kappa_2$ where we allow the empty cut which leaves the tree intact. Using $\tilde c$ we rewrite the last equation as
\begin{equation}
  \label{eq:ytau-step-2-4b}
  \begin{split}
 \der y^{[\tau^1 \cdots \tau^k]_a}_{ts} & =  \sum_{m \ge 0} \sum_{\zeta^1,\dots,\zeta^{k+m}} \sum_{\theta^1,\dots,\theta^m}
\frac{(k+m)!}{k!m!} \frac{\delta(\tau^1,\dots,\tau^k)}{\delta(\zeta^1,\dots,\zeta^{k+m})} \tilde c(\zeta^1,\tau^1,\theta^1) \cdots \tilde c(\zeta^k,\tau^k,\theta^k)
\\ & \qquad \times y^{\zeta} X^{\zeta^1 \cdots \zeta^{m} \theta^1 \cdots \theta^k}
   \end{split}
 \end{equation}
where now $\zeta = [\zeta^1\cdots \zeta^{k+m}]_a$ and $\zeta^1,\dots,\zeta^{k+m} \in \TT_\LL$ are non-empty trees and $\theta^1,\dots,\theta^k \in \FF_\LL$ are possibly empty forests but we exclude the case when $m=0$ and all the $\theta^i$ are empty. Now we will show that this expression corresponds exactly to
\begin{equation}
  \label{eq:ytau-step-2-5}
  \begin{split}
 \der y^{[\tau^1 \cdots \tau^k]_a}_{ts} & =  \sum_{m \ge 0} \sum_{\zeta \in \TT_\LL : \zeta = [\zeta^1\cdots \zeta^{k+m}]_a} c'(\zeta,[\tau^1 \cdots \tau^k]_a,\theta)  X^{\theta} y^\zeta
   \end{split}
 \end{equation}
which is what we want to prove. Note that the restriction in the sum over trees $\zeta$ of the form $[\zeta^1\cdots \zeta^{k+m}]_a$ for some $m \ge 0$ is due to the fact that for trees with less than $k$ branches at the origin the factor $c(\zeta,\tau,\theta)$ is zero.

 Each forest $\zeta^1\cdots \zeta^{k+m}$ appears
 $\delta(\zeta^1,\dots,\zeta^{k+m})$ times in the summation, moreover
 given the tree $\zeta = [\zeta^1\cdots \zeta^{k+m}]_a$ there are
 $(k+m)!/(k!m!)$ ways to choose $m$ branches of the root to cut
 away. Let us say that these cuts are on the last $m$ branches
 $\zeta^{k+1},\dots,\zeta^{k+m}$. Then the rest of the cuts appear on
 the first $k$ and for a fixed set $\zeta^1,\dots,\zeta^k$ of trees to
 cut there are $\delta(\tau^1,\dots,\tau^k)$ possible ways of
 associating each $\tau$ to some $\zeta$ to determine the associated
 cuts (if they are possible at all). Chosen the pairing between the
 $\zeta$-s and the $\tau$-s there are $\prod_{i=1}^k \sum_{\theta^i
   \in \FF_\LL} \tilde c(\zeta^i,\tau^i,\theta^i)$ possible cuts (note
 that chosen $\zeta^i$ and $\tau^i$ the forest $\theta^i$ is uniquely
 determined). Moreover since either $m> 0$ or some $\theta^i \neq
 \emptyset$ there is at least one proper (i.e. not empty nor full) cut
 in eq.~(\ref{eq:ytau-step-2-4b}). This concludes the proof.

\end{proof}

\section{Integration of finite increments}
\label{sec:one-dim}

We recall the integration theory introduced in~\cite{MR2091358} in
some details since this
setting is quite different from the original rough path
theory developed in~\cite{MR2036784,MR1654527}.

Notice that our future discussions will mainly rely on
$k$-increments with $k \le 3$. We measure the size of these increments by H\"older-like norms
: for $f \in \CC_2(V)$ let
$$
\norm{f}_{\mu} =
\sup_{s,t\in\ott}\frac{|f_{ts}|}{|t-s|^\mu},
\quad\mbox{and}\quad
\CC_1^\mu(V)=\lcl f \in \CC_2(V);\, \norm{f}_{\mu}<\infty  \rcl.
$$
In the same way, for $h \in \CC_3(V)$, set
\begin{eqnarray}
  \label{eq:normOCC2}
  \norm{h}_{\gamma,\rho} &=& \sup_{s,u,t\in\ott} 
\frac{|h_{tus}|}{|u-s|^\gamma |t-u|^\rho}\\
\|h\|_\mu &= &
\inf\left \{\sum_i \|h_i\|_{\rho_i,\mu-\rho_i} ;\, h =
 \sum_i h_i,\, 0 < \rho_i < \mu \right\} ,\nonumber
\end{eqnarray}
where the last infimum is taken over all sequences $\{h_i \in \CC_3(V) \}$ such that $h
= \sum_i h_i$ and for all choices of the numbers $\rho_i \in (0,z)$.
 We set
$$
\CC_3^\mu(V)=\lcl h\in\CC_3(V);\, \|h\|_\mu<\infty \rcl.
$$
Eventually,
let $\CC_3^{1+}(V) = \cup_{\mu > 1} \CC_3^\mu(V)$,
and remark that the same kind of norms can be considered on the
spaces $\cZ \CC_3(V)$, leading to the definition of the  spaces
$\cZ \CC_3^\mu(V)$ and $\cZ \CC_3^{1+}(V)$. 

\vspace{0.3cm}

With these notations in mind,
the following proposition is a basic result which is at the core of
 our approach to path-wise integration:
\begin{proposition}[The $\Lambda$-map]
\label{prop:Lambda}
There exists a unique linear map $\Lambda: \cZ \CC^{1+}_3(V)
\to \CC_2^{1+}(V)$ such that 
$$
\delta \Lambda  = \id_{\cZ \CC_3(V)}.
$$
Furthermore, for any $\mu > 1$, 
this map is continuous from $\cZ \CC^{\mu}_3(V)$
to $\CC_2^{\mu}(V)$ and we have
\begin{equation}\label{ineqla}
\|\Lambda h\|_{\mu} \le \frac{1}{2^\mu-2} \|h\|_{\mu} ,\qquad h \in  \cZ \CC^{1+}_3(V). 
\end{equation}
\end{proposition}

\vspace{0.3cm}

We can now give an algorithm for a canonical decomposition of the
preimage of $\cZ \CC_3^{1+}(V)$, or in other words, of a function 
$g\in\CC_2 (V)$ whose increment $\der g$ is small enough:
\begin{corollary}
\label{cor:integration}
Take an element $g\in\CC_2 (V)$ such that $\der g\in\CC_3^\mu(V)$ for
$\mu>1$. Then $g$ can be decomposed in a unique way as
$
g=\der f+ \Lambda \delta g,
$
where $f\in\CC_1(V)$. 

For any 2-increment $g\in\CC_2 (V)$, such that $\der g\in\CC_3^{1+}(V)$,
set
$
\delta f = (\id-\Lambda \delta) g 
$.
Then
$$
(\delta f)_{ts} = \lim_{|\Pi_{ts}| \to 0} \sum_{i=0}^n g_{t_{i+1}\, t_i},
$$
where the limit is over any partition $\Pi_{ts} = \{t_0=t,\dots,
t_n=s\}$ of $[t,s]$ whose size tends to zero.
\end{corollary}

\begin{proof}
See~\cite{MR2091358}.
\end{proof}

\section{Branched rough paths}
\label{sec:brp}

Up to this point we have considered only properties of the iterated
integrals of smooth functions $\{x^a\}_{a\in \LL}$ however from the
algebraic point of view the only data we need to build the family
$\{X^\tau\}_{\tau \in \TT_\LL}$ is a family of maps $\{I^a\}_{a \in
  \LL}$ from $\CC_2$ to $\CC_2$ satisfying certain properties.

\begin{definition}
\label{def:integral} We call \emph{integral} a linear map $I: \DD_I
   \to \DD_I$ on a sub-algebra $\DD_I \subset \CC^+_2$ satisfying two properties:
$$
I(hf) =  I(h) f, \qquad \forall h \in \DD_I, f \in \CC_1
$$
and 
$$
\der I(h) = I(e)h + \sum_i I(h^{1,i}) h^{2,i} \qquad \text{when $h\in
  \DD_I$, $ \der h = \sum_i h^{1,i} h^{2,i}$ and $h^{1,i} \in \DD_I$} 
$$
We explicitly require that $e \in \DD_I$.
\end{definition}

Using the embedding $f \in \CC_1 \mapsto fe \in \CC_2$ we can extend
the map $I$ to $\CC_1$: for any $f \in \CC_1$ we let $I(f) = I(fe)$
and since $fe = ef + \der f$  (as easily verified) we have
$$
I(f) = I(e)f + I(\der f) 
$$
for any $f \in \CC_1$ such that $\der f \in \DD_I$.

Given a family $\{I^a\}_{a\in\LL}$ of such integral maps on a common
algebra $\DD\subseteq \CC_{2}$ we can associate to them a family $\{X^{\tau}\}_{\tau\in\FF_\LL}$ recursively as done in Sect.~\ref{sec:itin} above:
$$
X^{\troot_a} = I^a(e), \qquad X^{[\tau^1 \cdots \tau^k]_a} = I^a(X^{\tau^1 \cdots \tau^k}), \qquad X^{\tau^1\cdots \tau^k}= X^{\tau^1}\circ \cdots \circ X^{\tau^k}.
$$   
In this way we estabilish an algebra homomorphism from $\cA\TT_\LL$ to
a subalgebra of $\CC_2$ generated by the $X^\tau$-s. This homomorphism
send the operation $B_+^a$ on $\cA\TT_\LL$ to the integral map $I^a$
on $\CC_2$. It is not difficult to verify that Theorem~\ref{th:equiv}
extends to the map $X$ generated by the family $\{I^a\}_a$.

Let us now introduce a regularity condition on the map $X$.

Given $\gamma \in (0,1]$ define the function $q_\gamma$ on forests as
 $q_\gamma(\tau) = 1$ for $|\tau| \le 1/\gamma$ and 
\begin{equation}
  \label{eq:q-def}
q_{\gamma}(\tau) = \frac{1}{2^{\gamma |\tau|}-2} \sum' q_\gamma(\tau^{(1)}) q_\gamma(\tau^{(2)})  
\end{equation}
whenever $\tau \in \TT$ with $|\tau| > 1/\gamma$ and  $q_\gamma(\tau_1 \cdots \tau_n) = q_\gamma(\tau_1) \cdots q_\gamma(\tau_n)$ for $\tau_1,\dots,\tau_n \in \TT$.

Note that $q_{\gamma}$ satisfy also the equation
$$
q_{\gamma}(\tau) = \frac{1}{2^{\gamma |\tau|}} \sum q_\gamma(\tau^{(1)}) q_\gamma(\tau^{(2)})
$$
which involves the splitting given by the coproduct $\Delta$ while the definition~(\ref{eq:q-def}) involves the splitting of trees given by the reduced coproduct $\Delta'$.

\begin{definition}
 We call an homomorphism $X :\cA \TT \to \CC_2$ a  \emph{branched
 rough path} (BRP) of roughness $\gamma> 0$, if it satisfy the equation~(\ref{eq:alg-relations})
 and moreover is such that
\begin{equation}
  \label{eq:brp-small}
 \| X^\tau\|_{\gamma |\tau|} \le B A^{|\tau|} q_\gamma(\tau), \qquad \tau \in \FF_\cl
\end{equation}
for some constants $B \in [0,1]$ and $A \ge 0$.
\end{definition}

Under certain conditions we can extend an homomorphism $X:
\cA_n \TT \to \CC_2$ defined only on the sub-algebra of trees with
degree less or equal to $n$ to the whole algebra.

\begin{theorem}
\label{th:ext}
Let us given a partial  homomorphism  $X:
\cA_n \TT_\cl \to \CC_2$ satisfying eq.~(\ref{eq:alg-relations}) and
such that there exists positive constants $\gamma,A\ge 0,B\in[0,1]$ for which 
\begin{equation}
  \label{eq:brp-small-2}
 \| X^\tau\|_{\gamma |\tau|} \le B A^{|\tau|} q_\gamma(\tau), \qquad \tau \in \TT^n_\cl
\end{equation}
with $\gamma (n+1) > 1$. Then there exists a unique extension of $X$
to a branched rough path defined on the whole $\cA \TT$ with roughness
$\gamma$ and such that eq.~(\ref{eq:brp-small-2}) holds for any $\tau
\in \TT_\LL$.
\end{theorem}

\begin{proof}
We  proceed by induction and assume that we have already found an extension $X:\cA_m \TT_{\cl} \to \CC_2$ satisfying eq.~(\ref{eq:alg-relations}) and for which we have
\begin{equation}
  \label{eq:brp-small-3}
 \| X^\tau\|_{\gamma |\tau|} \le B q_\gamma(\tau) A^{|\tau|}, \qquad \tau \in \TT^m_\cl.
\end{equation}
 This is true if $m=n$. Let us prove that we can extend $X$ to the set of trees with degree $m+1$ with the same bound on the H\"older norms. Since $\gamma m \ge \gamma (n+1) > 1$ we can set
$
X^\tau := \Lambda\left[ X^{\Delta' \tau}\right]   
$
for every $\tau$ such that $|\tau|=m$. Indeed
$$
\|X^{\Delta' \tau}\|_{m\gamma} \le  \sum_i \|X^{\tau^{(1)}_i \otimes \tau^{(2)}_i }\|_{m\gamma} \le  \sum_i'  \|X^{\tau^{(1)}_i}\|_{|\tau^{(1)}_i|\gamma} \|X^{ \tau^{(2)}_i }\|_{|\tau^{(2)}_i|\gamma}
$$
since $|\tau^{(1)}_i|+|\tau^{(2)}_i|=m$ for every $i$. This shows that $X^{\Delta' \tau} \in \CC_2^{m\gamma}$ and so it is in the domain of $\Lambda$.
To prove the bound on $X^\tau$ recall that
\begin{equation*}
  \begin{split}
\|X^\tau\|_{\gamma|\tau|}  & =   \|\Lambda X^{\Delta' \tau}\|_{\gamma|\tau|} \le \frac{1}{2^{|\tau|\gamma}-2}   
 \sum_i'  \|X^{\tau^{(1)}_i}\|_{|\tau^{(1)}_i|\gamma} \|X^{ \tau^{(2)}_i }\|_{|\tau^{(2)}_i|\gamma}
\\ & \le B^2 \frac{1}{2^{|\tau|\gamma}-2}   
 \sum_i' A^{|\tau^{(1)}_i|+|\tau^{(2)}_i|} q_\gamma(\tau^{(1)}_i) q_\gamma(\tau^{(2)}_i)
\\ & \le B^2 A^{|\tau|} q_\gamma(\tau)
  \end{split}
\end{equation*}
and since $B \le 1$ we have the required bound.
\end{proof}

\begin{remark}
\label{rem:q-asymp}
While we does not have been able to prove any asymptotic behavior for $q_\gamma(\tau)$ as $|\tau| \to \infty$ we conjecture that
\begin{equation}
  \label{eq:conject}
q_\gamma(\tau) \asymp C (\tau!)^{-\gamma}  
\end{equation}
for some constant $C$. For the class of linear Chen trees $\TT^{\text{Chen}}$ this conjecture is true thanks to the inequality
\begin{equation}
  \label{eq:neocx1}
\sum_{k=0}^n \frac{a^{\gamma k} b^{\gamma (n-k)}}{(k!)^\gamma (n!)^\gamma} \le c_\gamma \frac{(a+b)^{\gamma n}}{(n!)^\gamma}  
\end{equation}
valid for any $\gamma \in (0,1]$ and $a,b \ge 0$ and where the constant $c_\gamma$ depends only on $\gamma$. We prove this inequality in App.~\ref{app:neoc}. Note that this inequality is a variant of  Lyons' neo-classical inequality~(see e.g.\cite{MR1654527}) which in our notations reads
\begin{equation}
  \label{eq:neocxx1}
\sum_{k=0}^n \frac{a^{\gamma k} b^{\gamma (n-k)}}{(\gamma k)! [\gamma (n-k)n]!} \le c_\gamma \frac{(a+b)^{\gamma n}}{(\gamma n)!}  
\end{equation}
A sufficient condition for the validity of the conjecture would be the existence of a ``neo-classical tree inequality'' of the form
\begin{equation}
  \label{eq:neocx1-tree}
\sum \frac{a^{\gamma |\tau^{(1)}|} b^{\gamma |\tau^{(2)}|}}{(\tau^{(1)}!)^\gamma (\tau^{(2)}!)^\gamma} \le c_\gamma \frac{(a+b)^{\gamma |\tau|}}{(\tau !)^\gamma}  
\end{equation}
for any $\tau \in \TT$.
 The inequality is  true when $\gamma =1$ by using the tree binomial formula given in Lemma~\ref{lemma:tree-binomial}. 

The asymptotic behavior~(\ref{eq:conject}) appears also in the estimation of tree-indexed iterated integrals in the context of 3d Navier-Stokes equation studied in~\cite{MR2227041} (see also Sect.~\ref{sec:inf-dim}).
\end{remark}

We denote with $\Omega^\gamma_{\TT,\LL}$ the space of $\gamma$-BRP, on this
space we can introduce a distance by letting
$$
d_\gamma(X,Y) = \sum_{\tau \in \FF^n_\cl} \|X^\tau-Y^\tau\|_{\gamma|\tau|}
$$
where $n$ is again the largest integer such that $n\gamma \le 1$. This distance is strong enough to separate points in $\Omega^\gamma_{\TT,\LL}$: 

\begin{corollary}
If $X,Y \in \Omega^\gamma_{\TT,\LL}$ and $d_\gamma(X,Y) = 0$ then $X=Y$.  
\end{corollary}
\begin{proof}
If we let
$
Z^\tau = X^\tau-Y^\tau
$ for $\tau \in \TT^n_\LL$ then the partial homomorphism $Z$ is such that $Z^\tau = 0$ and satisfy
eq.~(\ref{eq:alg-relations}) for all $\tau \in \TT^n_\LL$. Then we
can choose $B=0$ and an arbitrary $A$ in  the
bounds~(\ref{eq:brp-small-2}) and use Thm.~\ref{th:ext} to conclude that we must have $Z^\tau =
0$ for any $\tau \in \TT_\LL$, i.e. that $X=Y$.  
\end{proof}

\begin{definition}
An~\emph{almost branched rough path} (aBRP) is a partial homomorphism $\widetilde X:
\cA_n \TT  \to \CC_2$ such that it approximately satisfy
eq.~(\ref{eq:alg-relations}) for any tree $\tau \in \TT_\cl^n$ modulus an element of $\CC^{1+}_3$
and for which we have
\begin{equation}
  \label{eq:abrp-small}
 \max_{\tau \in\TT_\cl^n} \|\widetilde X^\tau\|_{\gamma |\tau|} \le K 
\end{equation}
for some constant $K$ and some $\gamma > 1/(n+1)$.
 \end{definition}

Then we have the following result
\begin{theorem}
\label{th:corresp}
For any aBRP $\widetilde X$ there is a unique BRP $X$ of
roughness $\gamma$ such that 
$$
\max_{\tau \in \TT_\cl^n} \|X^\tau - \widetilde X^\tau\|_{(n+1)\gamma}
< \infty .
$$
\end{theorem}
\begin{proof}
The assumption is that $\der \widetilde X^\tau = \widetilde X^{\Delta' \tau} + R^\tau$ where $R^\tau \in \CC^{(n+1)\gamma}_3$ for any $\tau \in \FF_\cl^n$.

We will set $X^\tau = \widetilde X^\tau + Q^\tau$ and determine the increments $Q^\tau$ by induction. First look at $\tau$ such that $|\tau| = 1$, in this case
$$
\der X^\tau = \der \widetilde X^\tau + \der Q^\tau = R^\tau + \der Q^\tau
$$ 
since $\Delta' \tau = 0$. Then we set $Q^\tau = -\Lambda R^\tau$ since $R^\tau \in \ZZ\CC^{1+}_3$. So that we obtain $\der X^\tau = 0$ as it should.
Now assume that for $\tau \in \TT_\cl^{m}$ we have obtained $Q^\tau$ such that $\der X^\tau = X^{\Delta' \tau}$ and let us find such corrections $Q^{\tau}$ for $\tau \in \TT_\cl^{m+1}$ with $|\tau|=m+1$. We have
$$
\der \widetilde X^{\tau} = \sum' \widetilde X^{\tau^{(1)}} \widetilde X^{\tau^{(2)}} + R^\tau 
$$
since both $\tau^{(1)}$ and $\tau^{(2)}$ have degree less than $m+1$ we can apply the induction hypothesis and obtain
$
\der \widetilde X^{\tau} = \sum' ( X^{\tau^{(1)}}- Q^{\tau^{(1)}}) ( X^{\tau^{(2)}}- Q^{\tau^{(2)}}) + R^\tau.
$
Now let
$$
\widetilde R^\tau =
 \sum'  \left[Q^{\tau^{(1)}} X^{\tau^{(2)}}+X^{\tau^{(1)}}  Q^{\tau^{(2)}}-Q^{\tau^{(1)}}  Q^{\tau^{(2)}} \right]- R^\tau
$$
so that
$
\der \widetilde X^\tau - \widetilde R^\tau =  \sum'  X^{\tau^{(1)}} X^{\tau^{(2)}}.
$
If we can show that $\widetilde R^\tau \in \ZZ \CC_3^{1+}$, then setting $Q^{\tau} = \Lambda[\widetilde R^\tau]$ we would have obtained
$
\der X^\tau = \der \widetilde X^\tau - \widetilde R^\tau =  \sum'  X^{\tau^{(1)}} X^{\tau^{(2)}}.
$
as required and the induction would be complete.
It is clear that $\widetilde R^\tau \in \CC_3^{1+}$. The only problem is to prove that it is in the image of $\der$. By the triviality of the complex $(\CC_*,\der)$ this is equivalent to show that $\der \widetilde R^\tau = 0$. So let us prove the last equality.
Note that
$$
\der \widetilde R^\tau = \der\left[\der \widetilde X^\tau - \sum'  X^{\tau^{(1)}} X^{\tau^{(2)}}\right] =-  \der \sum'  X^{\tau^{(1)}} X^{\tau^{(2)}}
$$
Using again the induction hypothesis we get
$$
\der \widetilde R^\tau =   \sum'  X^{\tau^{(1)}} \der X^{\tau^{(2)}} -   \sum' \der X^{\tau^{(1)}} X^{\tau^{(2)}} 
$$
$$= \sum' X^{\tau^{(1)}} X^{\Delta' \tau^{(2)}} - \sum' X^{\Delta' \tau^{(1)}} X^{\tau^{(2)}} = X^{(\text{id}\otimes \Delta')\Delta' \tau} -X^{(\Delta'\otimes\text{id})\Delta' \tau} 
$$
But now
$
\der \widetilde R^\tau=  X^{(\text{id}\otimes \Delta')\Delta' \tau- (\Delta'\otimes\text{id})\Delta' \tau} = 0
$
since the reduced coproduct is coassociative. The proof of uniqueness is left to the reader.
\end{proof}

\section{Controlled paths}
\label{sec:controlled}
Following the line of development of~\cite{MR2091358} we describe now
a sufficiently large class of paths which can be integrated against a
given $\gamma$-branched rough path $X$. We then show that this set of
paths constitute an algebra and that integration and application of
sufficiently regular maps preserve this class. It will constitute the
natural space where to look for solutions of rough differential
equations driven by a branched path.  

In Sect.~\ref{sec:trees-diff} we showed that the solution  $y$ of a driven differential equation has the form of a series indexed by trees:
$
\der y_{ts} = \sum_{\tau \in \TT_\cl}  X^{\tau}_{ts} y^{\tau}_s
$ (cfr. eq.~(\ref{eq:tree-rep}))
 for suitable coefficients functions $\{ y^\tau : \tau \in \TT_\LL \}$ which satisfy eq.~(\ref{eq:ytau-der}). 

This suggest the following:

\begin{definition}
\label{def:controlled-path}
Let $X$ be a $\gamma$-BRP  and let $n$ the largest integer such that
$n\gamma\le 1$. For any $\kappa \in (1/(n+1),\gamma]$
a path $y$ is a $\kappa$-\emph{weakly controlled} by $X$ with values in the vector space $V$ if there exists paths $\{ y^{\tau} \in \CC^{|\tau|\kappa}_2(V) : \tau \in \FF_\LL^{n-1}\}$ and remainders $\{y^{\sharp} \in \CC_2^{n\kappa}(V), y^{\sharp,\tau} \in \CC_2^{(n-|\tau|)\kappa}(V), \tau \in \FF_{\cl}^{n-1}\}$ such that
\begin{equation}
  \label{eq:control}
\der y = \sum_{\tau \in \FF_\cl^{n-1}}  X^{\tau} y^{\tau} + y^\sharp  
\end{equation}
and for $\tau \in \FF_\cl^{n-1}$ :
\begin{equation}
  \label{eq:control2}
  \der y^{\tau} = \sum_{\sigma \in \FF_\cl^{n-1}}\sum_{\rho} c'(\sigma,\tau,\rho)  X^{\rho} y^{\sigma} + y^{\tau,\sharp}
\end{equation}
where we mean $\der y^{\tau} = y^{\tau,\sharp}$ when $|\tau| = n-1$.
We denote $\QQ_\kappa(X;V)$ the vector space  of $\kappa$-weakly
controlled paths by $X$ with values in $V$. Fixed a norm $|\cdot|$ on $V$   we introduce a norm $\|\cdot\|_{\QQ,\kappa}$ on $\QQ_\kappa(X;V)$ as
$$
\|y\|_{\QQ,\kappa} = |y_0| + \|y^\sharp\|_{n\kappa} + \sum_{\tau\in\FF_\LL^{n-1}} \|y^{\tau,\sharp}\|_{\kappa(n-|\tau|)}.
$$   
\end{definition}
To be precise, a well defined element in $\QQ_\kappa(X;V)$ is given by specifying the
path $y$ and all its ``derivatives'' $\{y^\tau\}_\tau$ but we usually
omit this for the sake of brevity.
A path in $\QQ_\kappa(X;V)$ has a partial
expansion in $X$ with a remainder denoted with $y^\sharp$. Likewise
every coefficient path in this expansion has a similar expansion of
progressively lower order.  We write $\QQ_{\kappa}(X) = \QQ_{\kappa}(X;\RR)$.

\begin{example}
\label{eq:example-2}  
Let us give an example with $d=1$ of the structure of a controlled path (since $d=1$ the partial series are indexed by unlabeled trees). Take $\gamma > 1/5$ so that $n=4$ and assume that
$X$ is a $\gamma$-BRP. Then $y \in \QQ_\gamma$ corresponds to the
set  of paths
\smalltrees
$$
y \in \CC_1^{\gamma},\qquad y^{\tsroot} \in \CC_1^{\gamma},\qquad 
y^{\aabb},y^{\tsnode \tsnode} \in \CC_1^{2\gamma}, \qquad    y^{\aababb},  y^{\aabb \tsnode}
y^{\aababb},   y^{\aaabbb}, y^{\tsnode \tsnode \tsnode}\in\CC_1^{3\gamma}
$$
which satisfy the following algebraic relations
\begin{equation*}
  \begin{split}
\der y & = X^{\tsroot} y^{\tsroot} + X^{\aabb} y^{\aabb} + X^{\tsnode \tsnode} y^{\tsnode \tsnode} +  X^{\aababb} y^{\aababb} + X^{\aabb \tsnode} y^{\aabb \tsnode} + X^{\aababb} y^{\aababb} + X^{\tsnode \tsnode \tsnode}  y^{\tsnode \tsnode \tsnode}+X^{\aaabbb} y^{\aaabbb} +y^\sharp  
\\
\der y^{\tsroot} & = X^{\tsroot} (y^{\aabb} + 2 y^{\tsroot \tsroot} ) + X^{\aabb} (y^{\aaabbb} + y^{\aabb \tsroot}) + X^{\tsroot \tsroot} (y^{\aabb \tsroot} + y^{\aababb} +3 y^{\tsroot \tsroot \tsroot}) + y^{\tsroot,\sharp}
\\
\der y^{\aabb} & = X^{\tsroot} (y^{\aabb \tsroot} + 2  y^{\aababb} + y^{\aaabbb}) + y^{\aabb,\sharp}
\\
\der y^{\tsroot \tsroot} & = X^{\tsroot} (y^{\aabb \tsroot}+y^{ \tsroot \tsroot \tsroot}) + y^{\tsroot \tsroot,\sharp}
\\
\der y^{\aababb} & =  y^{\aababb,\sharp} 
\\
\der y^{\aabb\tsnode} & =  y^{\aabb\tsnode,\sharp} 
\\
\der y^{\tsnode\tsnode\tsnode} & =  y^{\tsnode\tsnode\tsnode,\sharp} 
\\
\der y^{\aaabbb} & =  y^{\aaabbb,\sharp} 
  \end{split}
\end{equation*}
with remainders of orders
$$
y^{\sharp} \in \CC_2^{4\gamma},\quad y^{\tsroot,\sharp} \in
\CC_2^{3\gamma}, \qquad 
y^{\aabb,\sharp},y^{\tsroot \tsroot,\sharp} \in \CC_2^{2\gamma}
\qquad y^{\aababb,\sharp},y^{\aabb\tsnode,\sharp},y^{\tsnode\tsnode\tsnode,\sharp},y^{\aaabbb,\sharp}
\in \CC_2^{ \gamma } .
 $$
\end{example}
The following lemma will be useful in computations below.
\begin{lemma}
\label{lemma:derysharp}
\begin{equation*}
\der y^\sharp  = \sum_{\tau\in \FF_\cl^{n-1}}  X^\tau  y^{\tau,\sharp} 
\end{equation*}
\end{lemma}
\begin{proof}
  \begin{equation*}
  \begin{split}
\der y^\sharp & =  \sum_{\tau\in \FF_\cl^{n-1}}  X^\tau \der y^\tau - \sum_{\tau \in \FF_\cl^{n-1}}  \der X^\tau  y^\tau    
\\ & = \sum_{|\tau| = n-1}  X^\tau \der y^\tau+  \sum_{\tau\in \FF_\cl^{n-2}}  X^\tau \left(\sum_{\sigma\in \FF_\cl^{n-1}} \sum_\rho c'(\sigma,\tau,\rho)  X^{\rho} y^{\sigma} + y^{\tau,\sharp}\right) - \sum_{\sigma\in \FF_\cl^{n-1}}  \der X^\sigma  y^\sigma    
\\ & =  \sum_{|\tau| = n-1}  X^\tau \der y^\tau+ \sum_{\tau\in \FF_\cl^{n-2}}   \sum_{\sigma\in \FF_\cl^{n-1}} \sum_\rho c'(\sigma,\tau,\rho)  X^\tau X^{\rho} y^{\sigma} +  \sum_{\tau\in \FF_\cl^{n-2}}  X^\tau  y^{\tau,\sharp}  - \sum_{\sigma\in \FF_\cl^{n-1}}  \der X^\sigma  y^\sigma
\\ & =  \sum_{|\tau| = n-1}  X^\tau \der y^\tau+\sum_{\sigma\in \FF_\cl^{n-1}} \sum_{\tau \in \FF_\cl^{n-2},\rho} c'(\sigma,\tau,\rho)  X^\tau X^{\rho} y^{\sigma} +  \sum_{\tau\in \FF_\cl^{n-2}}  X^\tau  y^{\tau,\sharp}  - \sum_{\sigma\in \FF_\cl^{n-1}}  \der X^\sigma  y^\sigma
\\ & =\sum_{|\tau| = n-1}  X^\tau \der y^\tau+ \sum_{\sigma \in \FF_\cl^{n-1}} (X^{\Delta' \sigma} - \der X^\sigma) y^{\sigma} +  \sum_{\tau\in \FF_\cl^{n-2}}  X^\tau  y^{\tau,\sharp}  
\\ & = \sum_{\tau\in \FF_\cl^{n-1}}  X^\tau  y^{\tau,\sharp} 
  \end{split}
\end{equation*}
\end{proof}

\begin{lemma}
\label{lemma:controlled-phi}
Let $\varphi \in C_b^n(\RR^k,\RR)$ and $y \in \QQ_\kappa(X;\RR^k)$, then $z_t
= \varphi(y_t)$ is a weakly controlled path, $z \in \QQ_\kappa(X;\RR)$ where its coefficients are given by
\begin{equation*}
z^\tau = \sum_{m = 1}^{n-1} \sum_{\substack{\ol b \in \II\LL_1 \\ |\ol b|=m}}\frac{\varphi_{\ol b}(y)}{m!}
\sum_{\substack{\tau_1,\dots,\tau_m \in \FF_\cl^{n-1}\\ 
    \tau_1\cdots\tau_m=\tau }} y^{\tau_1,b_1} \cdots y^{\tau_m,b_m} , \qquad \tau \in \FF_\LL^{n-1}  
\end{equation*}
where $\LL_1=\{1,\dots,k\}$
(note that all the summations are over a finite number of terms).
\end{lemma}
\begin{proof}
The Taylor expansion for $\varphi$ reads
$$
\varphi(\xi') = \varphi(\xi) + \sum_{m = 1}^{n-1} \sum_{\substack{\ol b \in \II\LL_1 \\ |\ol b|=m}}\frac{\varphi_{\ol b}(\xi)}{m!} (\xi'-\xi)^{\ol b} + O(|\xi'-\xi|^n)
$$ 
which plugged into $\der z = \der \varphi(y)$ gives
\begin{equation*}
  \begin{split}
 \der z_{ts} & =  \sum_{m = 1}^{n-1} \sum_{\substack{\ol b \in \II\LL_1 \\ |\ol b|=m}}\frac{\varphi_{\ol b}(y_s)}{m!} (\der y_{ts})^{\ol b} + O(|t-s|^{n\kappa})
\\ & =  \sum_{m = 1}^{n-1} \sum_{\tau^1,\cdots\tau^m \in \FF_\LL^{n-1}} \sum_{\substack{\ol b \in \II\LL_1 \\ |\ol b|=m}}\frac{\varphi_{\ol b}(y_s)}{m!} y^{\tau^1,b_1}_s\cdots y^{\tau^m,b_m}_s X^{\tau^1\cdots \tau^m}_{ts} + O(|t-s|^{n\kappa})
  \end{split}
\end{equation*}
which gives the required result. To show that every $z^\tau$ satisfy the $\der$-equations~(\ref{eq:control2}) we can use a truncated version of the  arguments used in  Theorem~\ref{th:der-eqns-1}. We omit the details.
\end{proof}

The previous lemma shows that controlled paths are compatible with the
application of nonlinear functions. We will now prove that there
exists an extension of the  integral maps $\{I^a\}_a$ to the algebra
$\QQ_\gamma(X)$.

\begin{theorem}
\label{th:controlled-integrals}
The integral maps   $\{I^a\}_{a\in\LL}$ can be extended to  maps $I^a :  \QQ_\kappa(X) \to \der \QQ_\kappa(X)$.  If $y \in \QQ_\kappa(X)$ then $\der z = I^a( y)$ is such that
\begin{equation}
  \label{eq:2}
\der z = X^{\troot_a} z^{\troot_a} +    \sum_{\tau \in \TT_\cl^{n}} X^{\tau} z^{\tau} + z^{\flat}  
\end{equation}
where $z^{\troot_a} = y$,
$z^{[\tau]_a} = y^{\tau}$
and zero otherwise.
Moreover
$$
z^\flat =\Lambda\left[ \sum_{\tau  \in \FF_\cl^{n-1}\cup\{\emptyset\}}  X^{B^+_a(\tau)}  y^{\tau,\sharp}\right] \in \CC_2^{\kappa (n+1)}.
$$ 
\end{theorem}
\begin{proof}
Let $h = \sum_{\tau \in \FF^{n-1}_\LL} X^\tau y^\tau$ so that $\der y = h + y^\sharp$. By linearity and by the definition of $X$ we have $h \in \DD_I$ and 
$$
I^a(h) =\sum_{\tau \in \FF^{n-1}_\LL} I^a(X^{\tau}) y^\tau =  \sum_{\tau \in \FF^{n-1}_\LL} X^{[\tau]_a} y^\tau = I^a(\der y - y^\sharp)
$$  
we would like to show that we can extend $I^a$ such that $I^a(y^\sharp)$ is well defined so that we can set
$$
I^a(\der y) = \sum_{\tau \in \FF^{n-1}_\LL} X^{[\tau]_a} y^\tau + I^a(y^\sharp).
$$ 
To do this we compute the action of $\der$ on $I^a(y^\sharp)$. Since we want to preserve the properties of $I^a$ we have to require that
$$
\der I^a(y^\sharp) = I^a(e) y^\sharp + \sum_{\tau \in \FF^{n-1}_\LL}I^a(X^\tau) y^{\tau,\sharp}
= X^{\troot_a} y^\sharp + \sum_{\tau \in \FF^{n-1}_\LL}X^{[\tau]_a} y^{\tau,\sharp}
$$
where we used the computation of $\der y^{\tau,\sharp}$ in Lemma~\ref{lemma:derysharp}. Since $X$ is a $\gamma$-BRP and $y \in \QQ_\kappa(X)$ with $1/(n+1)< \kappa < \gamma$ we see that the r.h.s of this equation belongs to $\ZZ\CC_3^{(n+1)\kappa} \subset \ZZ\CC_3^{1+}$ so that it belongs to the domain of the $\Lambda$ map and then we can define
$$
I^a(y^\sharp) = \Lambda \left[X^{\troot_a} y^\sharp + \sum_{\tau \in \FF^{n-1}_\LL}X^{[\tau]_a} y^{\tau,\sharp}\right]
$$
which proves out statement taking into account that we can set $z = I^a(y) = I^a(1)y
+ I^a(\der y)$.
\end{proof}

\begin{example}
\label{ex:example-3}
Let us continue our one dimensional example. For the integral $z = I(y)$ of the controlled path $y$ introduced in Ex.~\ref{eq:example-2} we get
\smalltrees

\begin{equation*}
  \begin{split}
\der z = \der I(y) & = X^{\tsroot} y + X^{\aabb} y^{\tsroot} + X^{\aaabbb} y^{\aabb}+X^{\aababb} y^{\tsnode \tsnode} + X^{\aaabbabb} y^{\aabb \tsnode} + X^{\aaababbb} y^{\aababb} + X^{\aabababb}  y^{\tsnode \tsnode \tsnode}+X^{\aaaabbbb}y^{\aaabbb}+z^\flat  
\\ & =  X^{\tsroot} z^{\tsroot} + X^{\aabb} z^{\aabb} + X^{\aaabbb} z^{\aaabbb} +X^{\aababb} z^{\aababb} + z^{\sharp}
  \end{split}
\end{equation*}
with
$$
z^\flat = \Lambda\left[X^{\tsnode} y^{\sharp} + X^{\aabb} y^{\tsnode,\sharp} + X^{\aaabbb} y^{\aabb,\sharp} + X^{\aababb} y^{\tsnode \tsnode,\sharp} + X^{\aaababbb} y^{\aababb,\sharp}+X^{\aaabbabb} y^{\aabb \tsnode,\sharp} + X^{\aabababb} y^{\tsnode \tsnode \tsnode,\sharp}
\right].
$$
and the coefficients satisfy:
\begin{equation*}
  \begin{split}
\der z^{\tsroot} = \der y & = X^{\tsroot} y^{\tsroot} + X^{\aabb} y^{\aabb} + X^{\tsnode \tsnode} y^{\tsnode \tsnode} +  X^{\aababb} y^{\aababb} + X^{\aabb \tsnode} y^{\aabb \tsnode} + X^{\aababb} y^{\aababb} + X^{\tsnode \tsnode \tsnode}  y^{\tsnode \tsnode \tsnode}+X^{\aaabbb}y^{\aaabbb}+y^\sharp  
\\
& = X^{\tsroot} z^{\aabb} + X^{\aabb} y^{\aaabbb} + X^{\tsnode \tsnode} z^{\aababb} + z^{\tsroot,\sharp}
\\
\der z^{\aabb} = \der y^{\tsroot} & = X^{\tsroot} (y^{\aabb} + 2 y^{\tsroot \tsroot} ) + X^{\aabb} (y^{\aaabbb} + y^{\aabb \tsroot}) + X^{\tsroot \tsroot} (y^{\aabb \tsroot} + y^{\aababb} +3 y^{\tsroot \tsroot \tsroot}) + y^{\tsroot,\sharp}
\\
& = X^{\tsroot} (z^{\aaabbb} + 2 z^{\aababb}) + z^{\aabb,\sharp}
\\
\der z^{\aaabbb} = \der y^{\aabb} & = X^{\tsroot} (y^{\aabb \tsroot} + 2  y^{\aababb} + y^{\aaabbb}) + y^{\aabb,\sharp}
\\ & = z^{\aaabbb,\sharp}
\\
\der z^{\aababb} = \der y^{\tsroot \tsroot} & = X^{\tsroot} (y^{\aabb \tsroot}+y^{ \tsroot \tsroot \tsroot}) + y^{\tsroot \tsroot,\sharp}
\\ & = z^{\aababb,\sharp}
  \end{split}
\end{equation*}  
\end{example}

\begin{remark}
Given a controlled path $y \in \QQ_\kappa(X;\RR^n\otimes \RR^d)$ we can lift it to  a branched rough path $Y$ indexed by $\TT_{\LL_1}$ by the following recursion
$$
Y^{\troot_b} = \sum_{a \in \LL} I^a(y^{ab}), \qquad
Y^{[\tau^1 \cdots \tau^k]_b} = \sum_{a \in \LL} I^a(y^{ab} Y^{\tau^1}\circ \cdots \circ  Y^{\tau_k}), \qquad b \in \LL_1 
$$  
Indeed $\{ J^b(\cdot) = \sum_{a\in\LL}I^a(y^{ab} \cdot)\}_{b \in \LL_1}$ defines a family  of integrals in the sense of Def.~\ref{def:integral} and $Y$ is the associated $\gamma$-BRP.
\end{remark}

\subsection{Rough differential equations}

Let $f_a \in C(\RR^k;\RR^k)$, $a=1,\dots,d$ a family of vectorfields on $\RR^k$. Given a family on integral map $I^a$ which define a $\gamma$-BRP $X$ we consider the \emph{rough differential equation}
\begin{equation}
  \label{eq:rough-eq}
\der y = \sum_{a\in\LL} I^a(f_a(y)), \qquad y_0 = \eta \in \RR^k  
\end{equation}
in the time interval $[0,T]$. This equation has a well defined meaning
when the vectorfields $f_a$ are $C^n_b$ with $n$ the largest integer
for which $n\gamma \le 1$. In this case we can look for solutions
of the above equation with $y \in \QQ_\gamma(X;\RR^k)$ and
eq.~(\ref{eq:rough-eq}) can be understood as a fixed point problem in
$\QQ_\gamma(X;\RR^k)$ since we have that the map $\Gamma$ defined as
$$
\der \Gamma(y) = \sum_{a\in\LL} I^a(f_a(y)), \qquad \Gamma(y)_0 = \eta
$$ 
is well defined from $\QQ_\gamma(X;\RR^k)$ onto itself thanks to
Lemma~\ref{lemma:controlled-phi} and  Theorem~\ref{th:controlled-integrals}.

\begin{theorem}
If $\{f_a\}_{a\in\LL}$ is a family of $C^{n}_b$ vectorfields then
the rough differential equation~(\ref{eq:rough-eq}) has a global
solution $y \in \QQ_\gamma(X;\RR^k)$  for any initial condition
$\eta\in\RR^k$.

 If the vectorfields are
$C^{n+1}_b$ the solution 
$\Phi(\eta,X)\in\QQ_\gamma(X;\RR^k)$ is unique and the map $\Phi: \RR^k \times \Omega^\gamma_{\TT,\LL} \to
\QQ_{\gamma}(X;\RR^k)$ is Lipschitz in any finite interval $[0,T]$.   
\end{theorem}
\begin{proof}
The proof of existence is based on a compactness argument on the map
$\Gamma$. Global solutions are obtained exploiting the boundedness of
the vectorfields (and of their derivatives). 
Uniqueness is proven by contraction on sufficiently small
time interval $[0,S]$.  The arguments are just direct adaptation of
the proof of similar statements which can be found in~\cite{MR2091358}
and are quite standard so we prefer to omit them.
\end{proof}

\section{Infinite dimensional rough equations}
\label{sec:inf-dim}
Another motivation to introduce a rough path theory based on tree-indexed
iterated integrals comes from the observation that infinite
dimensional differential equations generate quite naturally expansions
in trees which cannot be reduced to ``linear'' iterated integrals by
the means of some geometric property. We still do not have a general
theory of such equations but in this section we would like to justify
our point of view by the means of three examples which we have studied
in detail elsewhere~\cite{kdv,MR2227041,TindelGubinelli}: the 1d periodic
deterministic Korteweg--de~Vries (KdV)  equation, Navier-Stokes like equations and a class of stochastic partial
differential equations. Given the illustrative purpose of this section
we will keep the exposition at a formal level. Rigorous results can be found in the papers cited above.

\subsection{The KdV equation}
The 1d periodic KdV equation is the partial differential equation
\begin{equation}
  \label{eq:kdv-real}
\partial_t u(t,\xi) + \partial^3_\xi u(t,\xi) +
\frac12 \partial_\xi u(t,\xi)^2 = 0,\quad u(0,\xi) = u_0(\xi), \qquad
(t,\xi)\in\RR\times \bTT  
\end{equation}
where  the initial condition $u_0$ belongs to some Sobolev space $H^\alpha(\bTT)$
of the torus $\bTT = [-\pi,\pi]$. This equation has many interesting
features (e.g. it is a completely integrable system) but here we are
interested only in the interplay between the non-linear term and the
dispersive linear term which is the generator of the Airy group
$U(t)$ of isometries
of $H^{\alpha}$. By going to Fourier variables and setting
$v_t = U(t) u_t$ we recast the above equation in integral
form
\begin{equation}
  \label{eq:kdv-base}
 v_t(k) =   v_0(k) + \frac{ik}2 \sum'_{k_1} \int_0^t e^{- i 3 k k_1 k_2 s}
  v_s(k_1) v_s(k_2) \, ds, \quad t \in [0,T], k \in \ZZ_* 
\end{equation}
where $k_2 = k-k_1$ and $v_0(k) = u_0(k)$ and where the primed
summation excludes the values $k_1 = 0$ and $k_1 = k$. We restrict our attention to
initial conditions such that $v_0(0) = 0$. By introducing the linear
operator
$
\dot X_\sigma(\varphi,\varphi) = \frac{ik}2 \sum'_{k_1} e^{- i 3 k k_1 k_2 \sigma}
  \varphi(k_1) \varphi(k_2)
$
this equation takes the abstract form
$$
v_t = v_s + \int_s^t \dot X_\sigma(v_\sigma,v_\sigma) d\sigma, \qquad
t,s \in [0,T].
$$
By iteratively substituting the unknown in this integral equation we
obtain an expansion whose first terms looks like
\begin{equation}
\label{eq:tree-series-kdv}
  \begin{split}
v_t & = v_s + \int_s^t d\sigma \dot X_\sigma(v_s,v_s)
+ 2 \int_s^t d\sigma \dot X_\sigma(v_s,\int_s^\sigma d\sigma_1 \dot
X_{\sigma_1}(v_s,v_s))
\\ &\qquad +\int_s^t d\sigma \dot X_\sigma(\int_s^\sigma d\sigma_1 \dot X_{\sigma_1}(v_s,v_s),\int_s^\sigma d\sigma_2 \dot X_{\sigma_2}(v_s,v_s))    
\\ &\qquad +4 \int_s^t d\sigma \dot X_\sigma(v_s,\int_s^\sigma
d\sigma_1 \dot X_{\sigma_1}(v_s,\int_s^{\sigma_1} d\sigma_2 \dot
X_{\sigma_2}(v_s,v_s))  + r_{ts}  
  \end{split}
\end{equation}
where $r_{ts}$ stands for the remaining terms in
the expansion. Denote with $\TT_{BP}\subseteq \TT$ the set of (unlabeled) planar rooted
trees with at most two branches at each node. A planar tree is a
rooted tree endowed with an ordering of the branches at each node. 
Then each of the terms in this expansion can be associated
to a tree in $\TT_{BP}$  and we can define recursively
multi-linear operators $X^{\tau}$ as
$$
X^{\troot}_{ts}(\varphi_1,\varphi_2) = \int_s^t \dot
X_\sigma(\varphi_1,\varphi_2) d\sigma;
$$ 
$$
X^{[\tau^1]}_{ts}(\varphi_1,\dots,\varphi_{m+1}) = \int_s^t \dot
X_\sigma(X^{\tau^1}_{\sigma
  s}(\varphi_1,\dots,\varphi_m),\varphi_{m+1})
 d\sigma
$$
and
$$
X^{[\tau^1 \tau^2]}_{ts}(\varphi_1,\dots,\varphi_{m+n}) = \int_s^t \dot
X_\sigma(X^{\tau^1}_{\sigma s}(\varphi_1,\dots,\varphi_m),X^{\tau^2}_{\sigma s}(\varphi_{m+1},\dots,\varphi_{m+n})) d\sigma.
$$
Eq.~(\ref{eq:tree-series-kdv}) has then the form
\smalltrees
\begin{equation}
  \label{eq:kdv-incr-1}
\der v_{ts} = X^{\tsroot}(v^{\times 2}) + X^{\aabb}(v^{\times 3}) +
X^{\aaabbb}(v^{\times 4}) +  X^{\aababb}(v^{\times 4}) +r  
\end{equation}
as an equation for $k$-increments where $v^{\times n}_s =
(v_s,\dots,v_s)$ ($n$ times). Moreover we have algebraic relations for
the $X^\tau$-s, for example
$$
\der X^{\aabb}(\varphi_1,\varphi_2,\varphi_3) = X^{\tsroot}(X^{\tsroot}(\varphi_1,\varphi_2),\varphi_3) ,
$$
$$
\der X^{\aaabbb}(\varphi_1,\varphi_2,\varphi_3,\varphi_4) =
X^{\tsroot}(X^{\aabb}(\varphi_1,\varphi_2,\varphi_3),\varphi_4) + X^{\aabb}(X^{\tsroot}(\varphi_1,\varphi_2),\varphi_3,\varphi_4) ,
$$
and
\begin{equation*}
  \begin{split}
\der X^{\aababb}& (\varphi_1,\varphi_2,\varphi_3,\varphi_4) =
X^{\tsroot}(X^{\tsroot}(\varphi_1,\varphi_2),X^{\tsroot}(\varphi_3,\varphi_4))
\\ & \qquad +
X^{\aabb}(\varphi_1,\varphi_2,X^{\tsroot}(\varphi_3,\varphi_4)) + 
X^{\aabb}(\varphi_3,\varphi_4,X^{\tsroot}(\varphi_1,\varphi_2))    
  \end{split}
\end{equation*}
where we used the symmetry of the operator $\dot X$ to obtain this
last equation. These relations have much in common with the analogous
relations for branched rough paths, however here the additional
information of the position of the various arguments must be taken
into account in the combinatorics of the reduced coproduct. It would
be interesting to determine  a Hopf algebra structure  on $\TT_{BP}$
which could account for these
relations in a general way. 

Our interest in the $X$-operators comes from the fact that they
are, usually,  more regular than the original operator $\dot X$. 
This additional regularity usually comes at the
expense of their H\oe lder time regularity when considered as 
operator-valued increments. We are then naturally led to consider
eq.~(\ref{eq:kdv-incr-1}) as a rough equation and to try to solve it
using the $\Lambda$ map. For example using only up to the double
iterated integrals we would obtain the equation
$$
\der v = (1-\Lambda \der)[X^{\tsroot}(v^{\times 2}) + X^{\aabb}(v^{\times 3})]
$$
which in some cases can be solved by fixed point methods.
This strategy has allowed us to obtain solutions of the KdV equation
for initial data in $H^{\alpha}$ with any $\alpha > -1/2$. Moreover it
provide a concrete strategy to improve this result in the sense that
if enough regularity of the two step-3 operators can be proven, then
we can solve the equation
$$
\der v = (1-\Lambda \der)[X^{\tsroot}(v^{\times 2}) +
X^{\aabb}(v^{\times 3}) +
X^{\aaabbb}(v^{\times 4}) +  X^{\aababb}(v^{\times 4}))]
$$
and  obtain solution for more irregular initial conditions.

\subsection{Navier-Stokes-like equations}

The $d$-dimensional NS equation (or the Burgers' equation) have the abstract form
\begin{equation}
  \label{eq:NS-c-abstract}
u_t = S_t u_0 +  \int_0^t S_{t-s}  B(u_s,u_s)\,ds.
\end{equation}
where $S$ is a bounded semi-group on a Banach space $\mathcal{B}$ and $B$ is a symmetric bilinear operator which is usually defined only on a subspace of $\mathcal{B}$. 
 Here we cannot proceed as in the
previous section since $S$ is only a semi-group and we must cope with
the convolution directly.
In~\cite{MR2227041} we showed that the solutions of this equation in the case of the 3d NS equation have
the series representation
\begin{equation}
  \label{eq:ns-series}
u_t = S_t u_0 + \sum_{\tau \in \TT_{B}} X^{\tau}_{t0}(u_0^{\times d(\tau)})   
\end{equation}
where $d(\tau)$ is a degree function and the $d(\tau)$-multilinear operator
$X^{\tau}$ has recursive definition
$$
X^{\tsroot}_{ts} (\varphi^{\times 2}) = \int_s^t S_{t-u}
B(S_{u-s}\varphi,S_{u-s}\varphi) du
$$
$$
X^{[\tau^1]}_{ts}(\varphi^{\times (d(\tau^1)+1)}) = \int_s^t
S_{t-u}B(X^{\tau^1}_{us}(\varphi^{\times d(\tau^1)}),\varphi) du
$$
and
$$
X^{[\tau^1 \tau^2]}_{ts}(\varphi^{\times (d(\tau^1)+d(\tau^2))}) = \int_s^t
S_{t-u}B(X^{\tau^1}_{us}(\varphi^{\times d(\tau^1)}),X^{\tau^2}_{us}(\varphi^{\times d(\tau^2)})) du
$$
These operators can be shown to allow bounds in $\mathcal{B}$ of the form
$$
|X^{\tau}(\varphi^{\times d(\tau)})|_{\mathcal{B}} \le C \frac{|t-s|^{\eps |\tau|}}{(\tau!)^{\eps}} |\varphi|_{\mathcal{B}}^{d(\tau)}
$$
where $\eps\ge 0$ is a constant depending on the particular Banach
space $\mathcal{B}$ we choose. The series~(\ref{eq:ns-series}) can be shown to be norm convergent at least for  small $t$ and define local solution of NS.
Due to the presence of the convolution integral these $X$ operators
does not behaves nicely with respect to
$\der$. In~\cite{TindelGubinelli} we introduced cochain complex $(\hat
C_*,\hat \der)$ adapted to the study of such convolution integrals
where the coboundary is given by $\tilde \der h = \der h - ah - ha$ with $a_{ts} =
S_{t-s}-\text{Id}$ the 2-increment naturally associated to the
semi-group. There exists also a corresponding $\tilde \Lambda$-map which
provide an appropriate inverse to $\tilde \der$. Algebraic relations for
these iterated integrals have then  by-now familiar expressions, e.g.:
$$
\tilde \der X^{\aabb}(\varphi^{\times 3}) =
X^{\tsroot}(X^{\tsroot}(\varphi^{\times 2}),\varphi) 
$$
etc...

\subsection{Polynomial SPDEs}
In the paper~\cite{TindelGubinelli} we study path-wise solutions to
SPDEs in the mild form
\begin{equation}
  \label{eq:SPDE}
u_t = S_t u_0 +  \int_0^t S_{t-s}  dw_s f(u_s)
\end{equation}
where the solution $u_t$ lives in some Hilbert space $\mathcal{B}$,
$S$ is an analytic semi-group in $\mathcal{B}$, $f:\mathcal{B} \to
\mathcal{V}$ some nonlinear function with values another Hilbert space
$\mathcal{V}$ and $w$ a Gaussian stochastic process with values in the
space of linear operators from $\mathcal{V}$ to $\mathcal{B}$
(possibly unbounded).
Like in the NS-like case above this abstract equation allows an
expansion in trees when the non-linear term is polynomial. For example
taking $f(\varphi) = B(\varphi,\varphi)$ for some symmetric bilinear
operator $B$ we get a stack of iterated integrals on the stochastic
process $w$:
$$
X^{\tsroot}_{ts} (\varphi^{\times 2}) = \int_s^t S_{t-u}
dw_u B(S_{u-s}\varphi,S_{u-s}\varphi) 
$$
$$
X^{[\tau^1]}_{ts}(\varphi^{\times (d(\tau^1)+1)}) = \int_s^t
S_{t-u} dw_u B(X^{\tau^1}_{us}(\varphi^{\times d(\tau^1)}),\varphi) 
$$ 
and
$$
X^{[\tau^1 \tau^2]}_{ts}(\varphi^{\times (d(\tau^1)+d(\tau^2))}) = \int_s^t
S_{t-u} dw_u B(X^{\tau^1}_{us}(\varphi^{\times d(\tau^1)}),X^{\tau^2}_{us}(\varphi^{\times d(\tau^2)})) 
$$
Where these integrals can be defined by stochastic integration with
respect to the process $w$ (It\^o or Stratonovich).
So provided useful (path-wise) estimates for these operators are available we can
use the $(\hat \CC,\tilde \der)$ complex and the $\tilde
\Lambda$ map  to set up rough equations and study path-wise solutions
of polynomial SPDE like eq.~(\ref{eq:SPDE}).



\appendix
\section{A variant of Lyons' neo-classical inequality}
\label{app:neoc}
\begin{proposition}
For any $\gamma \in (0,1]$ there exists a constant $c_\gamma$ such that
\begin{equation}
  \label{eq:neoc1}
\sum_{k=0}^n \frac{a^{\gamma k} b^{\gamma (n-k)}}{(k!)^\gamma ((n-k)!)^\gamma} \le c_\gamma \frac{(a+b)^{\gamma n}}{(n!)^\gamma}  
\end{equation}
for any $a,b > 0$.
\end{proposition}
\begin{proof}
Using Stirling's asymptotic for the factorial:
$
n! = e^{n(\log n -1)}\sqrt{2 \pi n} (1+O(1/n))
$  as $n \to \infty$ we can bound the sum $S_n$ on the l.h.s. of eq.~(\ref{eq:neoc1}) by
\begin{equation*}
S_n \le \frac{a^{\gamma n}}{(n!)^\gamma} + \frac{b^{\gamma n }}{(n!)^\gamma} +  \sum_{k=1}^{n-1} a^{\gamma k} b^{\gamma (n-k)} \frac{e^{\gamma k(1-\log k) + \gamma (n-k)(1-\log (n-k)) +d}}{(2\pi)^{\gamma} k^{\gamma} (n-k)^{\gamma}} g(k) 
\end{equation*}
where $g \ge 1$ is a bounded function such that $g(k) \to 1$ as $k \to \infty$ and $n-k \to \infty$.
Let $\varphi(x) = x \log(x/a) + (1-x) \log[(1-x)/b] + \log(a+b)$, then
$$
(n!)^\gamma (a+b)^{-\gamma n}  S_n \le \left(\frac{a}{a+b}\right)^{\gamma n} + \left(\frac{b}{a+b}\right)^{\gamma n } +  \sum_{k=1}^{n-1} (n!)^{\gamma} \frac{e^{\gamma (n - \log n) -\gamma n \varphi(k/n)}}
{(2\pi)^{\gamma} k^{\gamma} (n-k)^{\gamma}} 
$$
Using again the asymptotic formula for $n!$ we get
\begin{equation}
  \label{eq:asymp-1}
(n!)^\gamma (a+b)^{-\gamma n} S_n \le 2 +  \sum_{k=1}^{n-1} \frac{n^{\gamma} e^{-\gamma n \varphi(k/n)}}
{(2\pi)^{\gamma/2} k^{\gamma} (n-k)^{\gamma}} g'(k)   
\end{equation}
Where $g'$ is another function with the same properties as $g$. The function $\varphi$ has minimum in $a/(a+b)$ and $\varphi(a/(a+b)) = 0$. In the limit $n \to \infty$ the contributions to the sum coming from the terms for which $|k/n - a/(a+b)| > \eps$ is exponentially suppressed. Moreover $\varphi''(a/(a+b)) = (a+b)^2/(ab) \ge 1$ so the sum for the values of $k$ for which $|k/n - a/(a+b)| \le \eps$ can be bounded by a Gaussian integral uniformly in $a,b$. Then the r.h.s. of eq.~(\ref{eq:asymp-1}) can be bounded by a constant independent of $a,b$.
\end{proof}

\begin{remark}
The same approach can be used to prove the original neo-classical inequality if we do not care for  optimality of the  constant.
\end{remark}

\end{document}